\documentclass{svjour3}                     

\smartqed  

\usepackage{graphicx,psfrag}
\usepackage{amstext,amssymb,amsmath,amscd}
\usepackage{color,url}
\usepackage[curve]{xypic}
\usepackage{subfigure}
\usepackage{enumerate}

\newcommand{\R} {\mathbb R}

\newcommand{\p}{\partial}

\newcommand{\cancel}[1]{}
\newcommand{\raw}{\rightarrow}

\newcommand{\eps}{\epsilon}

\newcommand{\vn}[1]{\left|\left|#1\right|\right|}
\newcommand{\vsn}[1]{\left|#1\right|}

\newcommand{\ds}{\displaystyle}

\newcommand{\pent}{{\Omega_\epsilon}}
\newcommand{\pentp}{{\Omega'_P}}
\newcommand{\wwach}{w^{\rm Wach}}
\newcommand{\ltri}{\lambda^{\rm Tri}}
\newcommand{\lwach}{\lambda^{\rm Wach}}
\newcommand{\lopt}{\lambda^{\rm Har}}
\newcommand{\lsib}{\lambda^{\rm Sibs}}

\newcommand{\triang}{{\mathcal{T}}}

\newcommand{\Itri}{I^{\rm Tri}}
\newcommand{\Iwach}{I^{\rm Wach}}
\newcommand{\Iopt}{I^{\rm Opt}}
\newcommand{\Isib}{I^{\rm Sibs}}

\newcommand\diam{\textnormal{diam}}
\newcommand\dist{\textnormal{dist}}
\newcommand{\bc}{\textbf{c}}
\newcommand{\bp}{\textbf{p}}
\newcommand{\bq}{\textbf{q}}
\newcommand{\bv}{\textbf{v}}
\newcommand{\bx}{\textbf{x}}
\newcommand{\by}{\textbf{y}}
\newcommand{\bz}{\textbf{z}}
\newcommand{\bal}{{\bf \alpha}}

\newcommand{\lpn}[3]{\vn{#1}_{L^{#2}(#3)}}
\newcommand{\hpn}[3]{\vn{#1}_{H^{#2}(#3)}}
\newcommand{\hpsn}[3]{\vsn{#1}_{H^{#2}(#3)}}

\journalname{Advances in Computational Mathematics}

\begin{document}

\title{Error Estimates for Generalized Barycentric Interpolation\thanks{This research was supported in part by NIH contracts R01-EB00487, R01-GM074258, and a grant from the UT-Portugal CoLab project.}
}

\author{Andrew Gillette \and Alexander Rand \and Chandrajit Bajaj}

\authorrunning{Gillette, Rand, and Bajaj} 

\institute{A. Gillette \at
           Department of Mathematics\\
	   University of Texas at Austin\\
           \email{agillette@math.utexas.edu}
	   \and
	  A. Rand \at
	  Institute for Computational Engineering and Sciences \\
	  University of Texas at Austin\\
          \email{arand@ices.utexas.edu}
          \and
         C. Bajaj \at
         Department of Computer Science, Institute for Computational Engineering and Sciences \\
	 University of Texas at Austin\\
         \email{bajaj@cs.utexas.edu}
}

\date{Received: date / Accepted: date}



\thispagestyle{empty}
\begin{center}
{\LARGE\textbf{Error Estimates for Generalized Barycentric Interpolation}} \\
\textsc{Andrew Gillette\footnote{Department of Mathematics, University of Texas at Austin, \url{agillette@math.utexas.edu}}, Alexander Rand\footnote{Institute for Computational Engineering and Sciences, University of Texas at Austin, \url{arand@ices.utexas.edu}}, Chandrajit Bajaj\footnote{Department of Computer Science, Institute for Computational Engineering and Sciences, University of Texas at Austin, \url{bajaj@cs.utexas.edu}}} \\
\today
\end{center}

\begin{abstract}
We prove the optimal convergence estimate for first order interpolants used in finite element methods based on three major approaches for generalizing barycentric interpolation functions to convex planar polygonal domains.  The Wachspress approach explicitly constructs rational functions, the Sibson approach uses Voronoi diagrams on the vertices of the polygon to define the functions, and the Harmonic approach defines the functions as the solution of a PDE.  We show that given certain conditions on the geometry of the polygon, each of these constructions can obtain the optimal convergence estimate.  In particular, we show that the well-known maximum interior angle condition required for interpolants over triangles is still required for Wachspress functions but not for Sibson functions.

\end{abstract}

\section{Introduction}

While a rich theory of finite element error estimation exists for meshes made of triangular or quadrilateral elements, relatively little attention has been paid to meshes constructed from arbitrary polygonal elements.  Many quality-controlled domain meshing schemes could be simplified if polygonal elements were permitted for dealing with problematic areas of a mesh and finite element methods have been applied to such meshes~\cite{TS06,WBG07}.  Moreover, the theory of Discrete Exterior Calculus has identified the need for and potential usefulness of finite element methods using interpolation methods over polygonal domain meshes (e.g. Voronoi meshes associated to a Delaunay domain mesh)~\cite{GB2010}.  Therefore, we seek to develop error estimates for functions interpolated from data at the vertices of a polygon $\Omega$.

Techniques for interpolation over polygons focus on generalizing barycentric coordinates to arbitrary $n$-gons; this keeps the degrees of freedom associated to the vertices of the polygon which is exploited in nodal finite element methods.  The seminal work of Wachspress~\cite{W1975} explored this exact idea and has since spawned a field of research on rational finite element bases over polygons.  Many alternatives to these `Wachspress coordinates' have been defined as well, including the Harmonic and Sibson interpolants.  To our knowledge, however, no careful analysis has been made as to which, if any, of these interpolation functions provide the correct error estimates required for finite element schemes.

We consider first-order interpolation operators from some generalization of barycentric coordinates to arbitrary convex polygons.  A set of barycentric coordinates $\{\lambda_i\}$ for $\Omega$ associated with the interpolation operator $I: H^2(\Omega) \raw \text{span} \{\lambda_i\} \subset H^1(\Omega)$ is given by
\begin{equation}
\label{eq:genintop}
Iu:=\sum_iu(\bv_i)\lambda_i.  
\end{equation}
Since barycentric coordinates are unique on triangles (described in Section~\ref{ssec:triangulation}) this is merely the standard linear Lagrange interpolation operator when $\Omega$ is a triangle.

Before stating any error estimates, we fix some notation.  
For multi-index $\bal = (\alpha_1, \alpha_2)$ and point $\bx = (x,y)$, define $\bx^\bal := x^{\alpha_1} y^{\alpha_2}$, $\alpha ! := \alpha_1 \alpha_2$, $|\bal| := \alpha_1 + \alpha_2$, and $D^\bal u := \p^{|\bal|} u/\p x^{\alpha_1}\p y^{\alpha_2}$.  
The Sobolev semi-norms and norms over an open set $\Omega$ are defined by
\begin{align*}
\hpsn{u}{m}{\Omega}^2 &:=  \int_\Omega \sum_{|\alpha| = m} |D^\alpha u(\bx)|^2 \,{\rm d} \bx &{\rm and} & & \hpn{u}{m}{\Omega}^2 &:= \sum_{0\leq k\leq m}\hpsn{u}{m}{\Omega}^2.
\end{align*}
The $H^{0}$-norm is the $L^2$-norm and will be denoted $\lpn{\cdot}{2}{\Omega}$. 

Analysis of the finite element method often yields bounds on the solution error in terms of the best possible approximation in the finite-dimensional solution space.  
Thus the challenge of bounding the solution error is reduced to a problem of finding a good interpolant.  
In many cases Lagrange interpolation can provide a suitable estimate which is asymptotically optimal.  
For first-order interpolants that we consider, this \textbf{optimal convergence estimate} has the form
\begin{equation}
\label{eq:hconv}
\hpn{u - I u }{1}{\Omega} \leq C\,\diam(\Omega)\hpsn{u}{2}{\Omega},\quad\forall u\in H^2(\Omega).
\end{equation}

To prove estimate (\ref{eq:hconv}) in our setting, it is sufficient (see Section~\ref{sec:intsobsp}) to restrict the analysis to a class of domains with diameter one and show that $I$ is a bounded operator from $H^2(\Omega)$ into $H^1(\Omega)$, that is
\begin{equation}
\label{eq:iuh1uh2}
\hpn{Iu}{1}{\Omega} \leq C_I \hpn{u}{2}{\Omega},\quad \forall u\in H^1(\Omega).
\end{equation}
We call equation (\ref{eq:iuh1uh2}) the \textbf{$H^1$-interpolant estimate} associated to the barycentric coordinates $\lambda_i$ used to define $I$.

The optimal convergence estimate (\ref{eq:hconv}) does not hold uniformly over all possible domains; a suitable geometric restriction must be selected to produce a uniform bound.  Even in the simplest case (Lagrange interpolation on triangles), there is a gap between geometric criteria which are simple to analyze (e.g. the minimum angle condition) and those that encompass the largest possible set of domains (e.g. the maximum angle condition).  

\begin{figure}[ht]
\begin{center}
\psfrag{vi}{$\bv_i$}
\psfrag{Bi}{$\beta_i$}
\psfrag{c}{\textcolor{red}{$\bc$}}
\psfrag{diam}{\textcolor{blue}{$\diam(\Omega)$}}
\psfrag{Om}{$\Omega$}
\psfrag{r}{\textcolor{red}{$\rho(\Omega)$}}
\includegraphics[width=.3\linewidth]{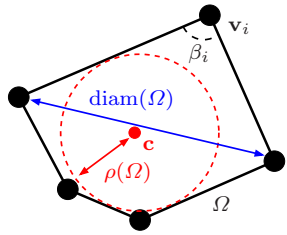}
\end{center} 
\caption{Notation used throughout paper. }
\label{fig:notation}
\end{figure}

This paper is devoted to finding geometric criteria under which the optimal convergence estimate (\ref{eq:hconv}) holds for several types of generalized barycentric coordinates on arbitrary convex polygons.  We begin by establishing some notation (shown in Figure \ref{fig:notation}) to describe the specific geometric criteria.  

Let $\Omega$ be a convex polygon with $n$ vertices.  Denote the vertices of $\Omega$ by $\bv_i$ and the interior angle at $\bv_i$ by $\beta_i$.  The largest distance between two points in $\Omega$ (the diameter of $\Omega$) is denoted $\diam(\Omega)$ and the radius of the largest inscribed circle is denoted $\rho(\Omega)$.  The center of this circle is denoted $\bc$ and is selected arbitrarily when no unique circle exists.  The \textbf{aspect ratio} (or chunkiness parameter) $\gamma$ is the ratio of the diameter to the radius of the largest inscribed circle, i.e.
\[\gamma := \frac{\diam(\Omega)}{\rho(\Omega)}.\]

We will consider domains satisfying one or more of the following geometric conditions.

\renewcommand{\labelenumi}{G\arabic{enumi}.}
\begin{enumerate}
\item \textbf{Bounded aspect ratio:} There exists $\gamma^*\in\R$ such that $\gamma < \gamma^*$.  \label{g:ratio}
\item \textbf{Minimum edge length: } There exists $d_*\in\R$ such that $|\bv_i - \bv_j| > d_* > 0$ for all $i\neq j$. \label{g:minedge}
\item \textbf{Maximum interior angle:} There exists $\beta^*\in\R$ such that $\beta_i < \beta^* < \pi$ for all $i$.\label{g:maxangle}
\end{enumerate}

Using several definitions of generalized barycentric functions from the literature, we show which geometric constraints on $\Omega$ are either necessary or sufficient to ensure the estimate for each definition.  The main results of this paper are summarized by the following theorem and Table \ref{tab:conditions}.  Primary attention is called to the difference between Wachspress and Sibson coordinates: while G\ref{g:maxangle} is a necessary requirement for Wachspress coordinates, it is demonstrated to be unnecessary for the Sibson coordinates.  

\begin{theorem}
In Table \ref{tab:conditions}, any necessary geometric criteria to achieve the $H^1$ interpolant estimate (\ref{eq:iuh1uh2}) are denoted by N.  The set of geometric criteria denoted by S in each row are sufficient to guarantee the $H^1$ interpolant estimate (\ref{eq:iuh1uh2}).

\begin{table}[ht]
\centering
\sbox{\strutbox}{\rule{0pt}{0pt}}           
\begin{tabular}[.75\textwidth]{@{\extracolsep{\fill}} ccccc}
 && \begin{tabular}{c} G\ref{g:ratio} \\ (aspect \\ ratio)\end{tabular} & \begin{tabular}{c} G\ref{g:minedge} \\ (min edge \\ length)\end{tabular} & 
\begin{tabular}{c} G\ref{g:maxangle} \\ (max interior \\ angle)\end{tabular} \\[2mm]
\hline\hline\\[2mm]
Triangulated & $\ds\ltri$ & - & - & S,N  \\[2mm]
\hline\\[2mm]
Harmonic & $\ds\lopt$  & S & - & - \\[2mm]
\hline\\[2mm]
Wachspress & $\ds\lwach$ & S & S & S,N \\[2mm]
\hline\\[2mm]
Sibson & $\ds\lsib$ & S & S & - \\[2mm]
\hline
\end{tabular}
\caption{`N' indicates a necessary geometric criterion for achieving the $H^1$ interpolant estimate (\ref{eq:iuh1uh2}).  The set of criteria denoted `S' in each row, taken together, are sufficient to ensure the $H^1$ interpolant estimate (\ref{eq:iuh1uh2}).}
\label{tab:conditions}
\end{table}
\end{theorem}
In Section \ref{sec:genbary}, we define the various types of generalized barycentric coordinates, compare their properties, and mention prior applications.
In Section \ref{sec:geomcond}, we review some general geometric results needed for subsequent proofs.
In Section \ref{sec:intsobsp}, we give the relevant background on interpolation theory for Sobolev spaces and state some classical results used to motivate our approach.
In Section \ref{sec:trimeth}, we show that the simplest technique of triangulating the polygon achieves the estimate if and only if G\ref{g:maxangle} holds.
In Section \ref{sec:optimal}, we show that if harmonic coordinates are used, G\ref{g:ratio} alone is sufficient.
In Section \ref{sec:wachpress}, we show that Wachspress coordinates require G\ref{g:maxangle} to achieve the estimate but all three criteria are sufficient.
In Section \ref{sec:sibson}, we show that Sibson coordinates achieve the estimate with G\ref{g:ratio} and G\ref{g:minedge} alone.
We discuss the implications of these results and future directions in Section \ref{sec:conc}.

\section{Generalized Barycentric Coordinate Types}
\label{sec:genbary}

Barycentric coordinates on general polygons are any set of functions satisfying certain key properties of the regular barycentric functions for triangles.

\begin{definition}\label{def:barcoor}
Functions $\lambda_i:\Omega\raw\R$, $i=1,\ldots, n$ are \textbf{barycentric coordinates} on $\Omega$ if they satisfy two properties. 
\renewcommand{\labelenumi}{B\arabic{enumi}.}

\begin{enumerate}
\item \textbf{Non-negative}:  $\lambda_i\geq 0$ on $\Omega$.\label{b:nonneg}
\item \textbf{Linear Completeness}: For any linear function $L:\Omega\raw\R$,
$\ds L=\sum_{i=1}^{n} L(\bv_i)\lambda_i$.\label{b:lincomp}
\end{enumerate}
\end{definition}

\begin{remark}
Property B\ref{b:lincomp} is the key requirement needed for our interpolation estimates.  It ensures that the interpolation operation preserves linear functions, i.e. $IL = L$.  
\end{remark}
We will restrict our attention to barycentric coordinates satisfying the following invariance property.  
Let $T:\R^2 \rightarrow \R^2$ be a composition of rotation, translation, and uniform scaling transformations and let $\{\lambda^T_i\}$ denote a set of barycentric coordinates on $T\Omega$.  

\renewcommand{\labelenumi}{B\arabic{enumi}.}
\begin{enumerate}
\setcounter{enumi}{2}
\item \textbf{Invariance:} $\ds\lambda_i(\bx)=\lambda_i^T(T(\bx))$.\label{b:invariance}
\end{enumerate}

This assumption will allow estimates over the class of convex sets with diameter one to be immediately extended to generic sizes since translation, rotation and uniform scaling operations can be easily passed through Sobolev norms (see Section~\ref{sec:intsobsp}).  At the expense of requiring uniform bounds over a class of diameter-one domains rather than a single reference element, complications associated with handling non-affine mappings between reference and physical elements are avoided \cite{ABF02}.  

A set of barycentric coordinates $\{\lambda_i\}$ also satisfies these additional familiar properties:
\renewcommand{\labelenumi}{B\arabic{enumi}.}
\begin{enumerate}
\setcounter{enumi}{3}
\item \textbf{Partition of unity:} $\ds\sum_{i=1}^{n}\lambda_i\equiv 1$. \label{b:partition}
\item \textbf{Linear precision:} $\ds\sum_{i=1}^{n}\bv_i\lambda_i(\bx)=\bx$. \label{b:linprec}
\item \textbf{Interpolation:} $\ds\lambda_i(\bv_j) = \delta_{ij}$. \label{b:interpolation}

\end{enumerate}
The precise relationship between these properties and those defining the barycentric coordinates is given in the following proposition.  
\begin{proposition} The properties B\ref{b:nonneg}-B\ref{b:interpolation} are related as follows:
\label{prop:Brelation}
\vspace{-.07in}
\begin{enumerate}[(i)]
\item B\ref{b:lincomp} $\Leftrightarrow$ (B\ref{b:partition} and B\ref{b:linprec})
\item (B\ref{b:nonneg} and B\ref{b:lincomp}) $\Rightarrow$ B\ref{b:interpolation}
\end{enumerate}
\end{proposition}

\begin{proof}
Given B\ref{b:lincomp}, setting $L\equiv 1$ implies B\ref{b:partition} and setting $L(\bx)=\bx$ yields B\ref{b:linprec}.  Conversely, assuming B\ref{b:partition} and B\ref{b:linprec}, let $L(x,y)=ax+by+c$ where $a,b,c\in\R$ are constants.  Let $\bv_i$ have coordinates $(\bv_i^x,\bv_i^y)$.  Then
\begin{align*}
\sum_{i=1}^{n}L(\bv_i)\lambda_i(x,y) & =\sum_{i=1}^{n} (a\bv_i^x+b\bv_i^y+c)\lambda_i(\bx)\\
 & =a\left(\sum_{i=1}^{n}\bv_i^x\lambda_i(\bx)\right)+b\left(\sum_{i=1}^{n}\bv_i^y\lambda_i(\bx)\right)+c\left(\sum_{i=1}^{n}\lambda_i(\bx)\right)\\
 & = ax+by+c = L(x,y).
\end{align*}
A proof that B\ref{b:nonneg} and B\ref{b:lincomp} imply B\ref{b:interpolation} can be found in \cite[Corollary 2.2]{FHK2006}.  \qed
\end{proof}
Thus, while other definitions of barycentric coordinates appear in the literature, requiring only properties B\ref{b:nonneg} and B\ref{b:lincomp} is a minimal definition still achieving all the desired properties.  

In the following subsections, we define common barycentric coordinate functions from the literature.  Additional comparisons of barycentric functions can be found in the survey papers of Cueto et al. \cite{CSCMCD2003} and Sukumar and Tabarraei \cite{ST2004}.

\subsection{Triangulation Coordinates}
\label{ssec:triangulation}

The simplest method for constructing barycentric coordinates on a polygon is to triangulate the polygon and use the standard barycentric coordinate functions of these triangles.  Interpolation properties of this scheme are well known from the standard analysis of the finite element method over triangular meshes, but this construction serves as an important point of comparison with the alternative barycentric coordinates discussed later.

Let $\triang$ be a triangulation of $\Omega$ formed by adding edges between the $\bv_j$ in some fashion.  Define
\[\ltri_{i,\triang}:\Omega\raw \R\]
to be the barycentric function associated to $\bv_i$ on triangles in $\triang$ containing $\bv_i$ and identically 0 otherwise.  Trivially, these functions define a set of barycentric coordinates on $\Omega$.

\begin{figure}[ht]
\begin{center}
\psfrag{vi}{$\bv_i$}
\psfrag{tmin}{$\triang_m$}
\psfrag{tmax}{$\triang_M$}
\[\begin{array}{ccc}
\includegraphics[width=.3\linewidth]{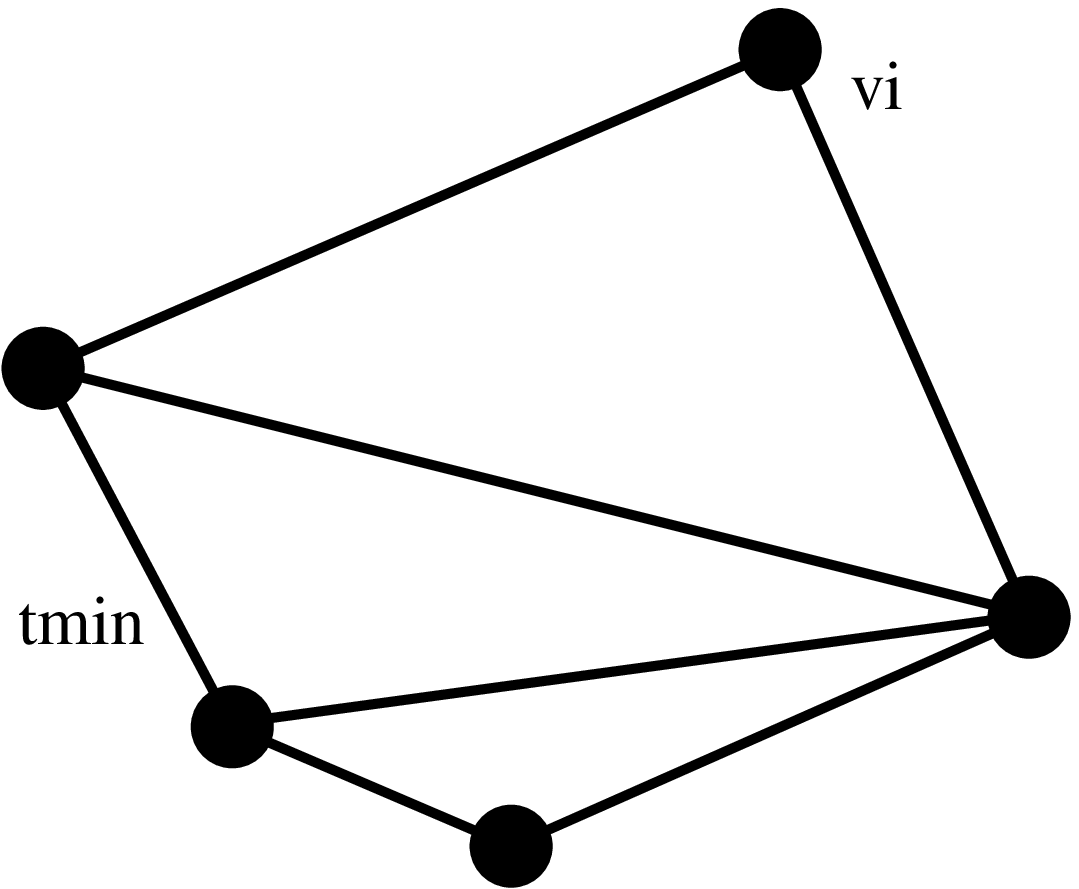} & 
\quad\text{  }\quad &
\includegraphics[width=.3\linewidth]{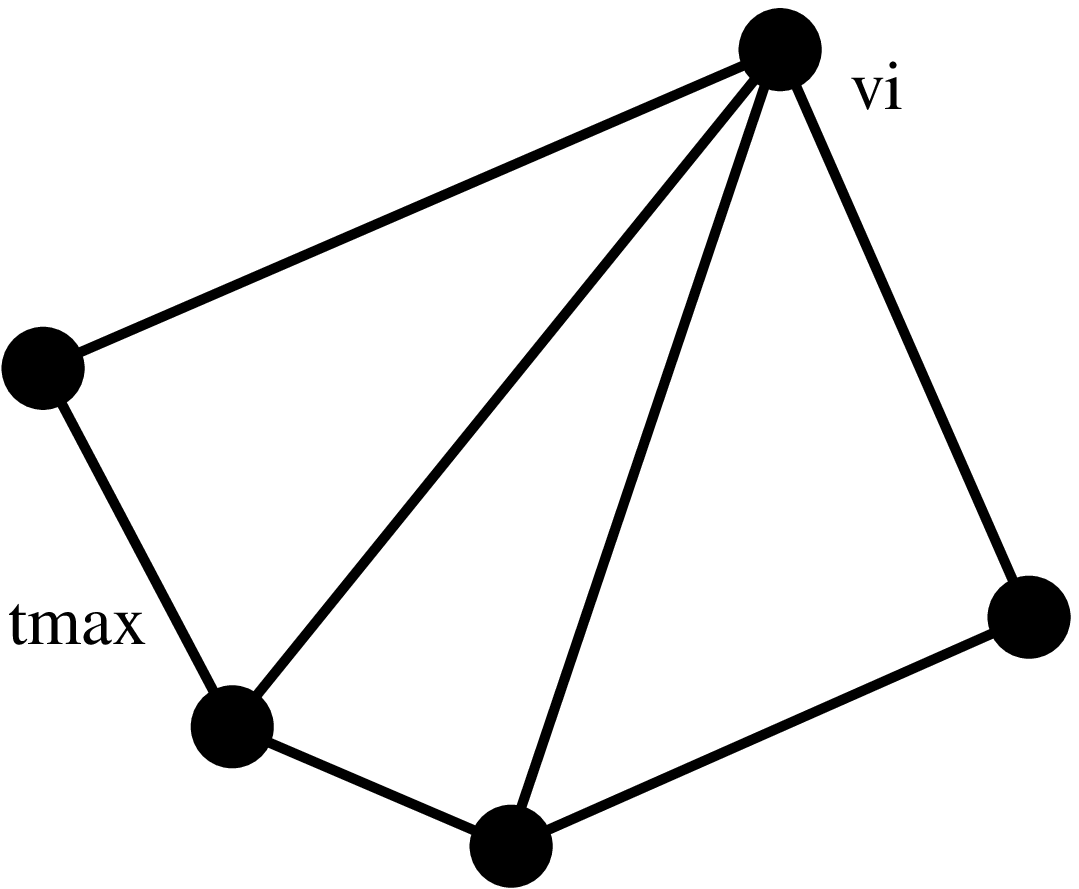}
\end{array}\]
\end{center} 
\caption{Triangulations $\triang_m$ and $\triang_M$ are used to produce the minimum and maximum barycentric functions associated with $\bv_i$, respectively.}
\label{fig:optimal}
\end{figure}

Two particular triangulations are of interest.  For a fixed $i$, let $\triang_m$ denote any triangulation with an edge between $\bv_{i-1}$ and $\bv_{i+1}$.  Let $\triang_M$ denote the triangulation formed by connecting $\bv_i$ to all the other $\bv_j$.  Examples are shown in Figure \ref{fig:optimal}. 

\begin{proposition}(Floater et al.~\cite{FHK2006}) 
\label{prop:fhk}
Any barycentric coordinate function $\lambda_i$ according to Definition \ref{def:barcoor} satisfies the bounds
\begin{equation}
\label{eq:fhk}
0\leq \ltri_{i,\triang_m}(\bx)\leq \lambda_i(\bx) \leq \ltri_{i,\triang_M}(\bx)\leq 1,\quad\forall\bx\in\Omega.
\end{equation}
\end{proposition}

Proposition \ref{prop:fhk} tells us that the triangulation coordinates are, in some sense, the extremal definitions of generalized barycentric coordinates.  In any triangulation of $\Omega$, at least one triangle will be of the form $(\bv_{i-1},\bv_i,\bv_{i+1})$, and hence the lower bound in (\ref{eq:fhk}) is always realized by some $\ltri_i$.  Thus, the examination of alternative barycentric coordinates can be motivated as an attempt to find non-extremal generalized barycentric coordinates.

\subsection{Harmonic Coordinates}

A particularly well-behaved set of barycentric coordinates, harmonic coordinates, can be defined as the solutions to certain boundary value problems.  Let $g_i:\p\Omega\raw\R$ be the piecewise linear function satisfying
\[g_i(\bv_j)=\delta_{ij},\quad g_i \text{ linear on each edge of $\Omega$}.\]
The harmonic coordinate function $\lopt_i$ is defined to be the solution of Laplace's equations with $g_i$ as boundary data,
\begin{equation}
\label{eq:optpde}
\ds\left\{\begin{array}{rcll}
\Delta\left(\lopt_i\right) & = & 0, & \text{on $\Omega$}, \\
\lopt_i & = & g_i. & \text{on $\p\Omega$}.
\end{array}\right.
\end{equation}
Existence and uniqueness of the solution are well known results~\cite{Ev98,RR04}.  Properties B\ref{b:nonneg} and B\ref{b:lincomp} are a consequence of the maximum principle and linearity of Laplace's equation.

These coordinates are optimal in the sense that they minimize the norm of the gradient over all functions satisfying the boundary conditions,
\[
\lopt_i = \text{argmin} \left\{\hpsn{\lambda}{1}{\Omega} \, : \, \lambda = g_i \,\text{on $\p\Omega$} \right\}.  
\]
This natural construction extends nicely to polytopes, as well as to a similar definition for barycentric-like (Whitney) vector elements on polygons.  Christensen~\cite{C2008} has explored theoretical results along these lines.  Numerical approximations of the $\lopt_i$ functions have been used to solve Maxwell's equations over polyhedral grids~\cite{E2007} and for finite element simulations for computer graphics~\cite{MKBWG2008,JMRGS07}.  While it may seem excessive to solve a PDE just to derive the basis functions for a larger PDE solver, the relatively limited geometric requirements required for their use (see Section \ref{sec:optimal}) make these functions a useful reference point for comparison with simpler constructions and a suitable choice in contexts where mesh element quality is hard to control.

\subsection{Wachspress Coordinates}

\begin{figure}[ht]
\begin{center}
\psfrag{vi}{$\bv_i$}
\psfrag{vip1}{$\bv_{i+1}$}
\psfrag{vim1}{$\bv_{i-1}$}
\psfrag{x}{$\bx$}
\psfrag{Aix}{$A_i(\bx)$}
\psfrag{Aip1x}{$A_{i+1}(\bx)$}
\psfrag{Aim1x}{$A_{i-1}(\bx)$}
\psfrag{wBi}{$B_i$}
\psfrag{dots}{$\cdots$}
\[\begin{array}{ccc}
\includegraphics[width=.3\linewidth]{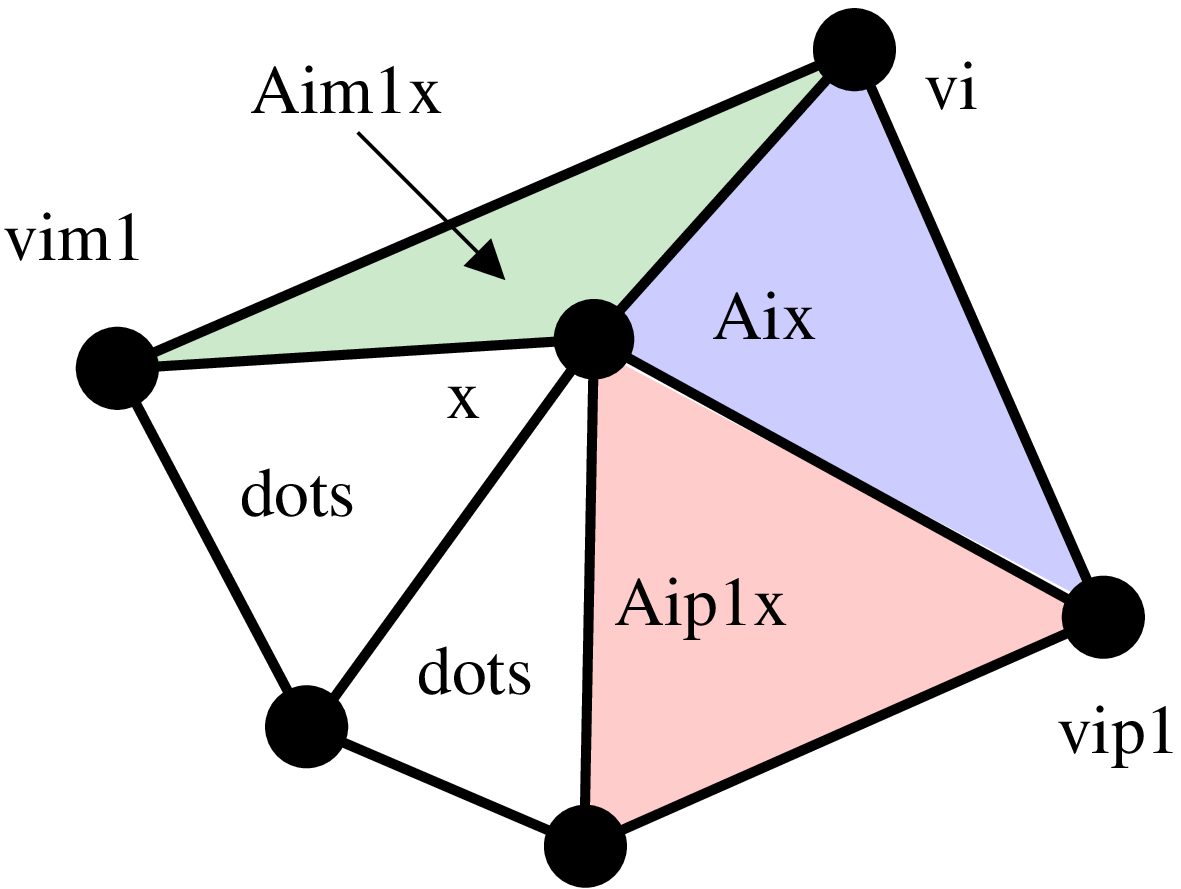} &
\includegraphics[width=.3\linewidth]{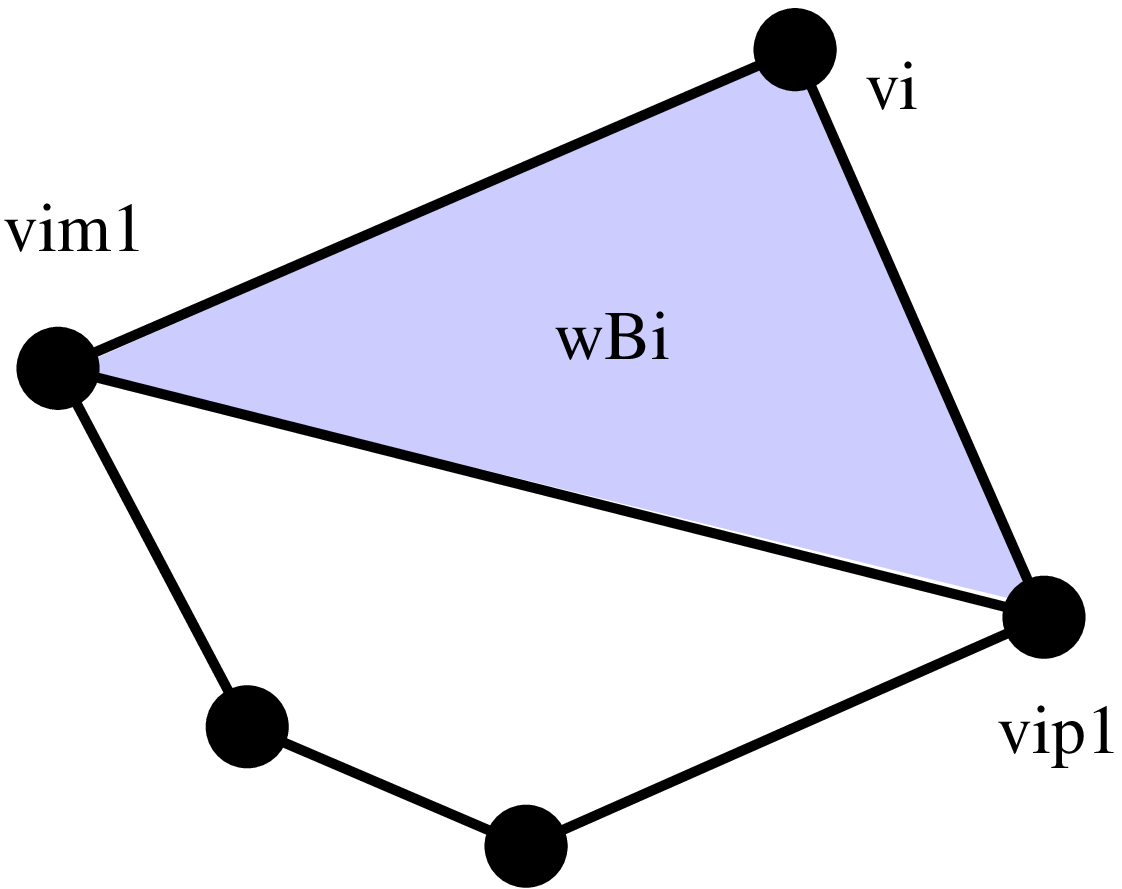}
\end{array}\]
\end{center} 
\caption{Left: Notation for $A_i(\bx)$. Right: Notation for $B_i$.}
\label{fig:wach-ntn}
\end{figure}

One of the earliest generalizations of barycentric coordinates was provided by Wachspress \cite{W1975}.  Definition of these coordinates is defined based on some notation shown in Figure \ref{fig:wach-ntn}.  Let $\bx$ denote an interior point of $\Omega$ and let $A_i(\bx)$ denote the area of the triangle with vertices $\bx$, $\bv_i$, and $\bv_{i+1}$ where, by convention, $\bv_{0} := \bv_{n}$ and $\bv_{n+1} := \bv_{1}$.  Let $B_i$ denote the area of the triangle with vertices $\bv_{i-1}$, $\bv_i$, and $\bv_{i+1}$. Define the Wachspress weight function
\[\wwach_i(\bx) = B_i \prod_{j\not=i,i-1}A_j(\bx).\]
The Wachspress coordinates are then given by
\begin{equation}
\label{eq:wach}
\lwach_i(\bx)=\frac{\wwach_i(\bx)}{\sum_{j=1}^{n} \wwach_j(\bx)}
\end{equation}
These coordinates have received extensive attention in the literature since they can be represented as rational functions in Cartesian coordinates.  Their use in finite element schemes has been numerically tested in specific application contexts but to our knowledge has not been evaluated in the general Sobolev error estimate context considered here.  We note that $\lwach_i\in H^1(\Omega)$ since it is a rational function with strictly positive denominator on $\Omega$.

%
%

\begin{remark}
Since $B_i$ does not depend on $\bx$ and $A_i(\bx)$ is linear in $\bx$, the Wachspress functions are degree $n-2$.  By a result from Warren~\cite{W2003}, the Wachspress functions are the unique, lowest degree rational barycentric functions over polygons.  For finite element applications, however, the $\lambda_i$ need not be rational.
\end{remark}

\subsection{Sibson (Natural Neighbor) Coordinates}

\begin{figure}[ht]
\psfrag{vi}{$\bv_i$}
\psfrag{Ci}{$C_i$}
\psfrag{x}{$\bx$}
\psfrag{Dix}{$D(\bx)$}
\psfrag{DixCi}{$D(\bx)\cap C_i$}
\[\begin{array}{ccc}
\includegraphics[width=.3\linewidth]{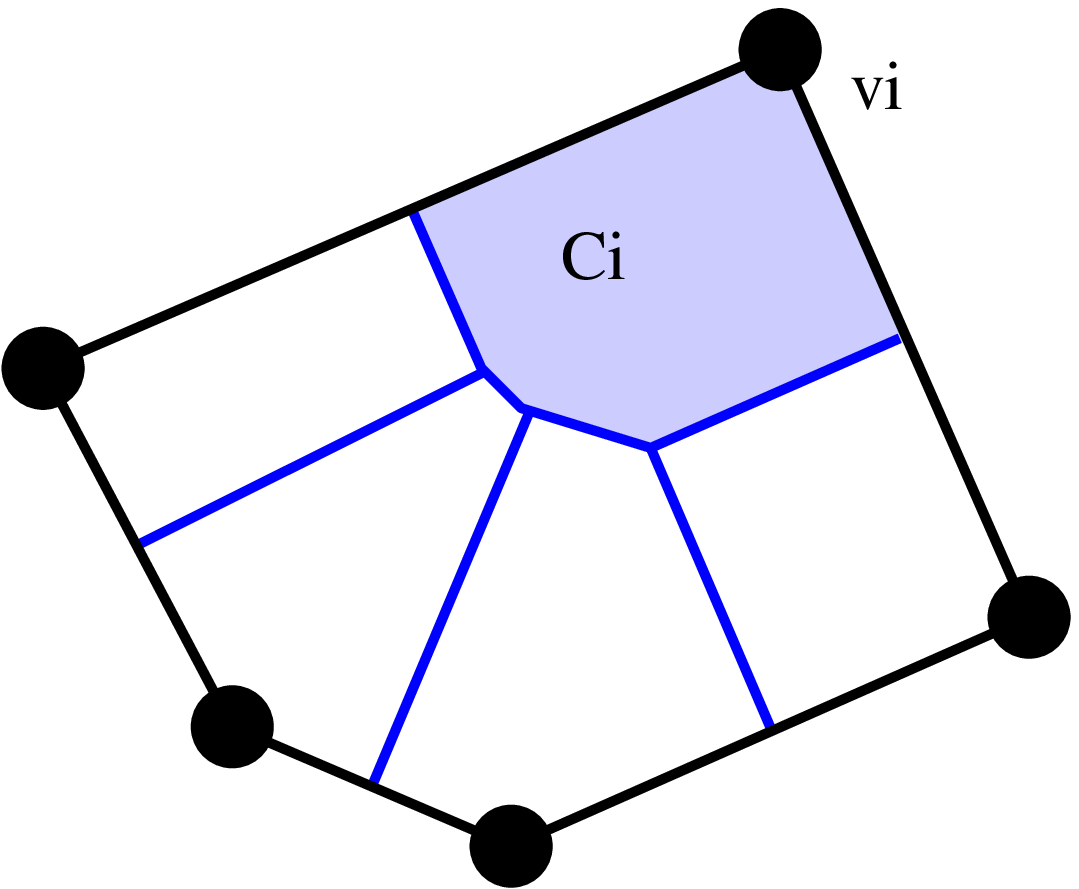} &
\includegraphics[width=.3\linewidth]{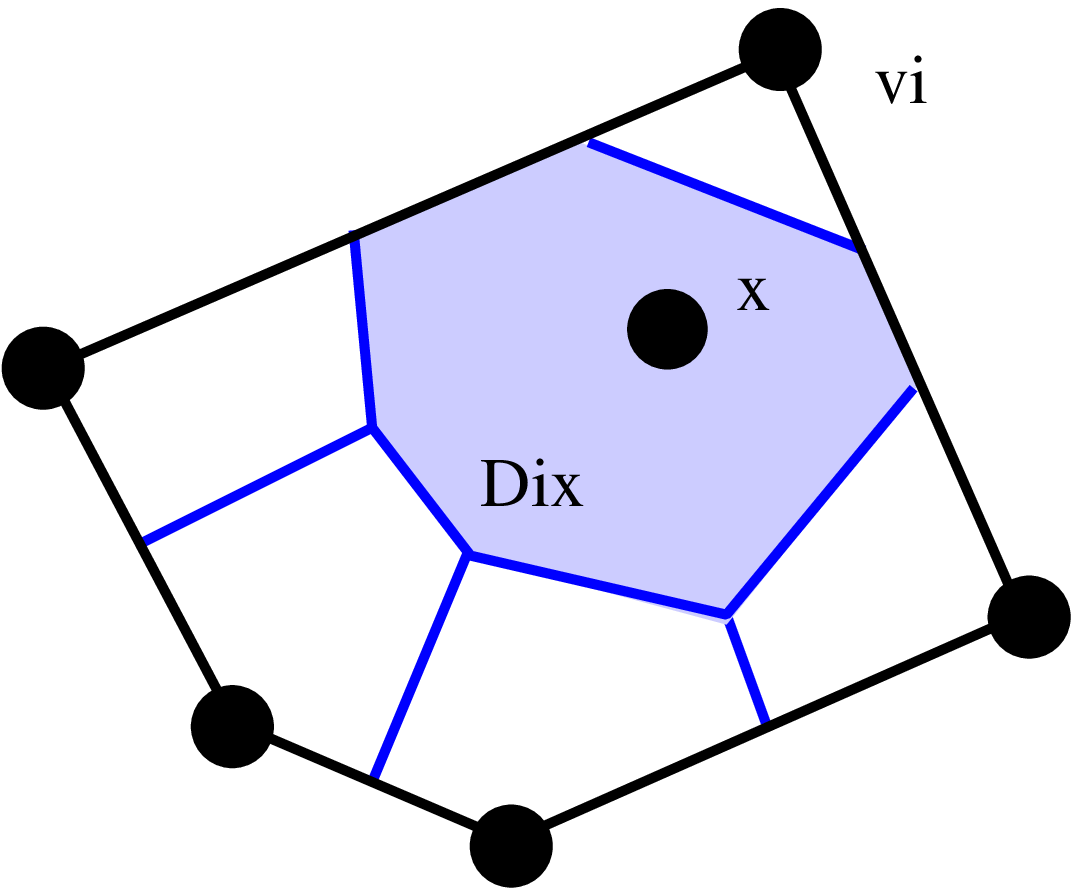} &
\includegraphics[width=.3\linewidth]{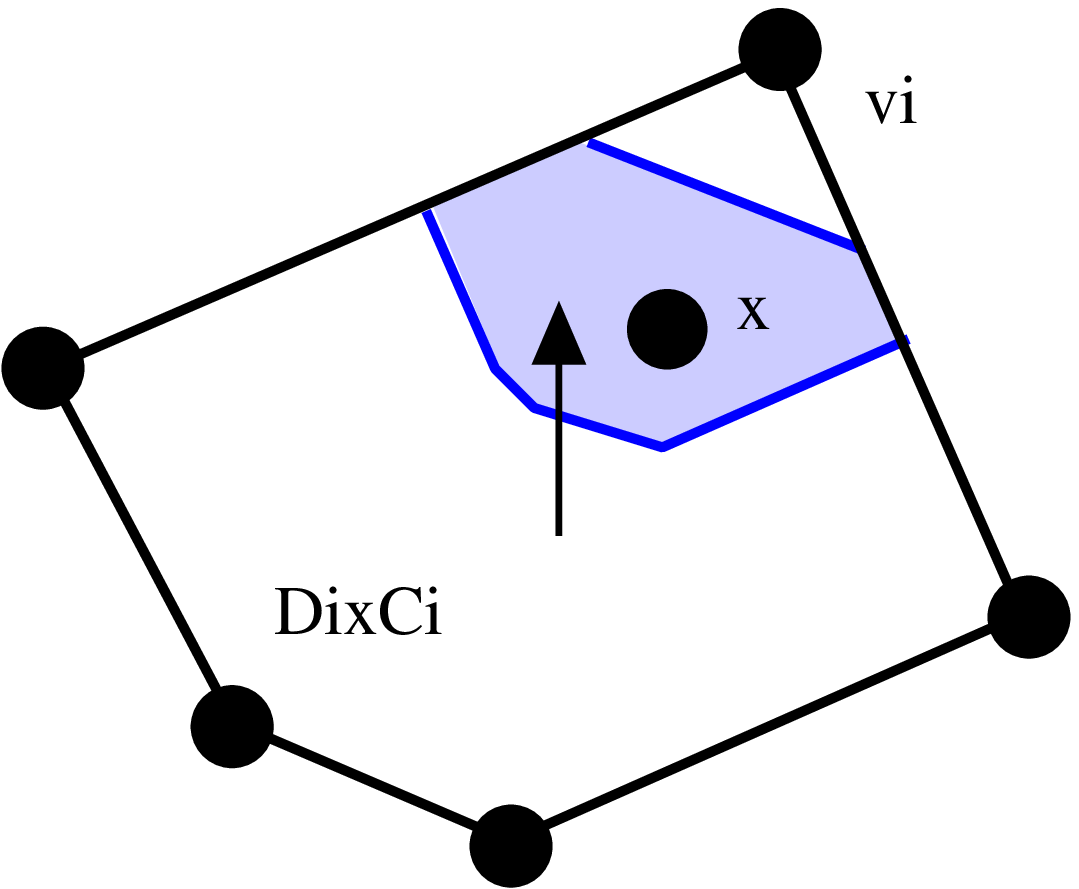}
\end{array}\]
\caption{Geometric calculation of a Sibson coordinate.  $C_i$ is the area of the Voronoi region associated to vertex $\bv_i$ inside $\Omega$.  $D(\bx)$ is the area of the Voronoi region associated to $\bx$ if it is added to the vertex list.  The quantity $D(\bx)\cap C_i$ is exactly $D(\bx)$ if $\bx=\bv_i$ and decays to zero as $\bx$ moves away from $\bv_i$, with value identically zero at all vertices besides $\bv_i$.}
\label{fig:sibson}
\end{figure}

The Sibson coordinates \cite{S1980}, also called the natural neighbor or natural element coordinates, make use of Voronoi diagrams on the vertices $\bv_i$ of $\Omega$.  Let $\bx$ be a point inside $\Omega$.  Let $P$ denote the set of vertices $\{\bv_i\}$ and define
\[P'=P\cup \{\bx\} = \{\bv_1,\ldots,\bv_{n},\bx\}.\]
We denote the \textbf{Voronoi cell} associated to a point $\bp$ in a pointset $Q$ by 
\[V_Q(\bp):= \left\{\by\in \Omega \, : \, \vsn{\by-\bp} < \vsn{\by-\bq} \, , \, \forall\bq\in Q\setminus\{\bp\}  \right\}.\]  
Note that these Voronoi cells have been restricted to $\Omega$ and are thus always of finite size.  We fix the notation
\[\begin{array}{lcccl}
C_i &:=& |V_{P}(\bv_i)| &=& \left|\{\by\in\Omega \, : \, \vsn{\by-\bv_i} < \vsn{\by-\bv_j}\, , \, \forall j\not=i\}\right| \\
 &&  &=& \text{area of cell for $\bv_i$ in Voronoi diagram on the points of $P$,} \\
\\
D(\bx) &:=&  |V_{P'}(\bx)| &=& \left|\{\by\in\Omega \, : \, \vsn{\by-\bx} < \vsn{\by-\bv_i}\, , \, \forall i\}\right| \\
 &&  &=& \text{area of cell for $\bx$ in Voronoi diagram on the points of $P'$}.
\end{array}
\]
By a slight abuse of notation, we also define
\[ D(\bx)\cap C_i := |V_{P'}(\bx)\cap V_{P}(\bv_i)|. \]
The notation is shown in Figure \ref{fig:sibson}.  The Sibson coordinates are defined to be
\begin{align*}
\lsib_i(\bx) &:= \frac{D(\bx)\cap C_i}{D(\bx)} &
\textnormal{or, equivalently,} & & 
\lsib_i(\bx) &= \frac{D(\bx)\cap C_i}{\sum_{j=1}^{n}D_j(\bx)\cap C_j}.
\end{align*}

It has been shown that the $\lsib_i$ are $C^\infty$ on $\Omega$ except at the vertices $\bv_i$ where they are $C^0$ and on circumcircles of Delaunay triangles where they are $C^1$ \cite{S1980,F1990}.  Since the finite set of vertices are the only points at which the function is not $C^1$, we conclude that $\lsib_i\in H^1(\Omega)$.  

To close this section, we compare the intra-element smoothness properties of the coordinate types on the interior of $\Omega$.  The triangulation coordinates are $C^0$, the Sibson coordinates are $C^1$, and the Wachspress functions and the harmonic coordinates are both $C^\infty$.

\section{Generalized Shape Regularity Conditions}
\label{sec:geomcond}

The invariance property B\ref{b:invariance} allows estimates on diameter-one polygons to be scaled to polygons of arbitrary size.  Several well-known properties of planar convex sets to be used throughout the analysis are given in Proposition~\ref{prop:convexfacts}.  Let $|\Omega|$ denote the area of convex polygon $\Omega$ and let $|\partial \Omega|$ denote the perimeter of $\Omega$.

\begin{proposition}\label{prop:convexfacts}
If $\Omega$ is a convex polygon with $\diam(\Omega) = 1$, then
\vspace{-.07in}
\begin{enumerate}[(i)]
\item $|\Omega| < \pi/4$,\label{cf:maxarea}
\item $|\partial \Omega| \leq \pi$, \label{cf:maxperimeter}
\item $\Omega$ is contained in a ball of radius no larger than $1/\sqrt{2}$, and \label{cf:jungs}
\item If convex polygon $\Upsilon$ is contained in $\Omega$, then $|\partial \Upsilon| \leq |\partial \Omega|$. \label{cf:subset}
\end{enumerate}
\end{proposition}

The first three statements are the isodiametric inequality, a corollary to Barbier's theorem, and Jung's theorem, respectively.  The last statement is a technical result along the same lines.  See \cite{Eg58,YB61,SA00} for more details.

Certain combinations of the geometric restrictions (G1-G3) imply additional useful properties for the analysis.  These resulting conditions are listed below.  

\renewcommand{\labelenumi}{G\arabic{enumi}.}
\begin{enumerate}
\setcounter{enumi}{3}
\item \textbf{Minimum interior angle:} There exists $\beta_*\in\R$ such that $\beta_i > \beta_* > 0$ for all $i$.\label{g:minangle}
\item \textbf{Maximum vertex count:} There exists $n^*\in\R$ such that $n < n^*$.\label{g:maxdegree}  
\end{enumerate}

For triangles, G\ref{g:minangle} and G\ref{g:maxangle} are the only two important geometric restrictions since G\ref{g:maxdegree} holds trivially and G\ref{g:ratio}$\Leftrightarrow$G\ref{g:minangle}$\Rightarrow$G\ref{g:minedge}.  For general polygons, the relationships between these conditions are more complicated; for example, a polygon satisfying G\ref{g:ratio} may have vertices which are arbitrarily close to each other and thus might not satisfy G\ref{g:maxdegree}.  Proposition~\ref{prop:Grelation} below specifies when the original geometric assumptions (G1-G3) imply G\ref{g:minangle} or G\ref{g:maxdegree}.

\begin{proposition} The following implications hold.
\label{prop:Grelation}
\vspace{-.07in}
\begin{enumerate}[(i)]
\item G\ref{g:ratio} $\Rightarrow$ G\ref{g:minangle} \label{gr:getminangle}
\item (G\ref{g:minedge} or G\ref{g:maxangle}) $\Rightarrow$ G\ref{g:maxdegree}\label{gr:getmaxdegree}
\end{enumerate}
\end{proposition}
\begin{proof}
G\ref{g:ratio} $\Rightarrow$ G\ref{g:minangle}: If $\beta_i$ is an interior angle, then $\rho(\Omega) \leq \sin(\beta_i/2)$ (see Figure~\ref{fig:minIntAngProof}).  Thus $\gamma > \frac{1}{\sin(\beta_i/2)}$.  We conclude that $\beta_i > 2 \arcsin \frac{1}{\gamma^*}$.  Note that $\gamma^*\geq 2$ so this is well-defined.

\begin{figure}[ht]
\begin{center}
\psfrag{vi}{$\bv_i$}
\psfrag{Bi}{$\beta_i$}
\psfrag{c}{\textcolor{red}{$\bc$}}
\psfrag{r}{\textcolor{red}{$\rho(\Omega)$}}
\includegraphics[width=.3\linewidth]{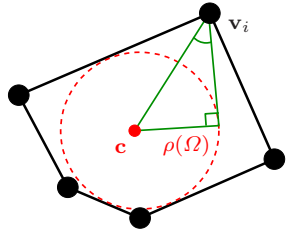}
\end{center} 
\caption{Proof that G1 $\Rightarrow$ G4.  The upper angle in the triangle is $\leq\beta_i/2\leq \pi/2$ and the hypoteneuse is $\leq\diam(\Omega)=1$.  Thus $\rho(\Omega)\leq \sin(\beta_i/2)$.}
\label{fig:minIntAngProof}
\end{figure}

G\ref{g:minedge} $\Rightarrow$ G\ref{g:maxdegree}: By Jung's theorem (Proposition~\ref{prop:convexfacts}(\ref{cf:jungs})), there exists $\bx\in\Omega$ such that $\Omega \subset B(\bx,1/\sqrt{2})$.  By G\ref{g:minedge}, $\{B(\bv_i,d_*/2)\}_{i=1}^{n}$ is a set of disjoint balls.  Thus $B(\bx, 1/\sqrt{2}+d_*/2)$ contains all of these balls.  Comparing the areas of $\bigcup_{i=1}^{n} B(\bv_i,d_*/2)$ and $B(\bx, 1/\sqrt{2}+d_*/2)$ gives $n\frac{\pi d_*^2}{4} < \pi (\frac 1{\sqrt 2}+d_*/2)^2$, so $n < \frac{(\sqrt{2}+d_*)^2}{d_*^2}$.

G\ref{g:maxangle} $\Rightarrow$ G\ref{g:maxdegree}: Since $\Omega$ is convex, $\sum_{i=1}^{n} \beta_i = \pi(n-2)$.  So $n \beta^* \geq \pi(n-2)$.  Thus $n \leq \frac{2\pi}{\pi-\beta^*}$. 
\qed
\end{proof}


\section{Interpolation in Sobolev Spaces}
\label{sec:intsobsp}



Interpolation error estimates are typically derived from the Bramble-Hilbert lemma which says that Sobolev functions over a certain domain or class of domains can be approximated well by polynomials.  
The original lemma \cite{BH70} applied to a fixed domain (typically the ``reference'' element) and did not indicate how the estimate was impacted by domain geometry.  Later, a constructive proof based on the averaged Taylor polynomial gave a uniform estimate under the geometric restriction G\ref{g:ratio} \cite{DS80,BS08}.  
Recent improvements to this construction have demonstrated that even the condition G\ref{g:ratio} is unnecessary \cite{Ve99,DL04}.  This modern version of the Bramble-Hilbert lemma is stated below and has been specialized to our setting, namely, the $H^1$ estimate for diameter $1$, convex domains.

\begin{lemma}[\cite{Ve99,DL04}]\label{lem:bramblehilbert}
Let $\Omega$ be a convex polygon with diameter $1$.  For all $u\in H^2(\Omega)$, there exists a first order polynomial $p_u$ such that 
$\hpn{u-p_u}{1}{\Omega} \leq C_{BH} \, \hpsn{u}{2}{\Omega}$. 
\end{lemma}

We emphasize that the constant $C_{BH}$ is uniform over all convex sets of diameter $1$.  The $H^1$-interpolant estimate (\ref{eq:iuh1uh2}) and Lemma~\ref{lem:bramblehilbert} together ensure the desired 
optimal convergence estimate (\ref{eq:hconv}).

\begin{theorem}
Let $\Omega$ be a convex polygon with diameter $1$.  If the $H^1$-interpolant estimate (\ref{eq:iuh1uh2}) holds, then for all $u\in H^2(\Omega)$, 
\[
\hpn{u-Iu}{1}{\Omega} \leq (1 + C_I)\, \sqrt{1+C_{BH}^2} \, \hpsn{u}{2}{\Omega}.
\]  
\end{theorem}
\begin{proof}
Let $p_u$ be the polynomial given in Lemma~\ref{lem:bramblehilbert} which closely approximates $u$.  By property B\ref{b:lincomp}, $Ip_u = p_u$ yielding the estimate
\begin{align*}
\hpn{u-Iu}{1}{\Omega}  & \leq \hpn{u-p_u}{1}{\Omega} + \hpn{I(u-p_u)}{1}{\Omega}\\
& \leq (1 + C_I)\hpn{u-p_u}{2}{\Omega}
 \leq (1 + C_I)\sqrt{1+C_{BH}^2} \hpsn{u}{2}{\Omega}.\qed
\end{align*}
\end{proof}
 
\begin{corollary} 
Let $\diam(\Omega) \leq 1$.  If the $H^1$-interpolant estimate (\ref{eq:iuh1uh2}) holds, then for all $u\in H^2(\Omega)$,  
\[ 
\hpn{u-Iu}{1}{\Omega} \leq (1 + C_I)\, \sqrt{1+C_{BH}^2}\, \diam(\Omega) \, \hpsn{u}{2}{\Omega}. 
\]   
\end{corollary} 
\begin{proof} 
This follows from the standard scaling properties of Sobolev norms since property B\ref{b:invariance} allows for a change of variables to a unit diameter domain.   
Note: the $L^2$-component of the $H^1$-norm satisfies a stronger estimate containing an extra power of $\diam(\Omega)$.   
\qed 
\end{proof}

Section~\ref{sec:errorest} is an investigation of the geometric conditions under which the $H^1$-interpolant estimate (\ref{eq:iuh1uh2}) holds for the barycentric functions discussed in Section~\ref{sec:genbary}.  
Under the geometric restrictions G\ref{g:ratio} and G\ref{g:maxdegree}, one method for verifying (\ref{eq:iuh1uh2}) (utilized in \cite{BS08} for simplicial interpolation) is to bound the $H^1$-norm of the basis functions.  In several cases we will utilize this criteria which is justified by the following lemma.  

\begin{lemma}
\label{lem:basisbd}
Under G\ref{g:ratio} and G\ref{g:maxdegree}, the $H^1$-interpolant estimate (\ref{eq:iuh1uh2}) holds whenever there exists a constant $C_\lambda$ such that
\begin{equation}\label{eq:basisbound}
\hpn{\lambda_i}{1}{\Omega} \leq C_\lambda.
\end{equation}
\end{lemma}
\begin{proof}
This follows almost immediately from the Sobolev embedding theorem; see \cite{Ad03,Le09}: 
\[
\hpn{Iu}{1}{\Omega} \leq \sum_{i=1}^{n} |u(\bv_i)| \hpn{\lambda_i}{1}{\Omega} \leq n^* C_\lambda\, \vn{u}_{C^0(\overline{\Omega})} \leq n^* \,C_\lambda\, C_s\, \hpn{u}{2}{\Omega},
\]
where $C_s$ is the Sobolev embedding; i.e., $\vn{u}_{C^0(\overline{\Omega})} \leq C_s\, \hpn{u}{2}{\Omega}$ for all $u\in H^2(\Omega)$.  The constant $C_s$ is independent of the domain $\Omega$ since the boundaries of all polygons satisfying G\ref{g:ratio} are uniformly Lipschitz~\cite{Le09}. \qed
\end{proof}


\section{Error Estimate Requirements}\label{sec:errorest}

\subsection{Estimate Requirements for Triangulation Coordinates}
\label{sec:trimeth}

Interpolation error estimates on triangles are well understood: the optimal convergence estimate (\ref{eq:hconv}) holds as long as the triangle satisfies a maximum angle condition \cite{BA76,Ja76}.  In fact, it has been shown that the triangle circumradius controls the error independent of any other geometric criteria \cite{Kr91}.  This result can be directly applied to $\Itri$, the interpolation operator associated to coordinates $\ltri_i$.  This convention will also be used to define $\Iopt$, $\Iwach$, and $\Isib$ as the interpolation operators associated with with harmonic, Wachspress, and Sibson coordinates, respectively.  

\begin{lemma}
\label{lem:triangulated}
Under G\ref{g:maxangle}, the $H^1$ interpolant estimate (\ref{eq:iuh1uh2}) holds for $\Itri$.  Conversely, G\ref{g:maxangle} is a necessary assumption to achieve (\ref{eq:iuh1uh2}) with $\Itri$.
\end{lemma}
\begin{proof}
All angles of all triangles of any triangulation $\triang$ of $\Omega$ satisfying G\ref{g:maxangle} are less than $\beta^*$.  Thus, the sufficiency of G\ref{g:maxangle} follows immediately from the maximum angle condition on simplices~\cite{BA76}.  An example from the same paper involving the interpolation of a quadratic function over a triangle also establishes the necessity of the condition.\qed
\end{proof}


\subsection{Estimate Requirements for Harmonic Coordinates}
\label{sec:optimal}



\begin{figure}[ht]
\begin{center}
\psfrag{vi}{$\bv_i$}
\psfrag{vip1}{$\bv_{i+1}$}
\psfrag{vim1}{$\bv_{i-1}$}
\psfrag{c}{\textcolor{red}{$\bc$}}
\psfrag{y}{\textcolor{red}{$\by$}}
\psfrag{x}{$\bx$}
\psfrag{Aix}{$A_i(\bx)$}
\psfrag{wBi}{$B_i$}
\[\begin{array}{cccc}
\includegraphics[height=.2\linewidth]{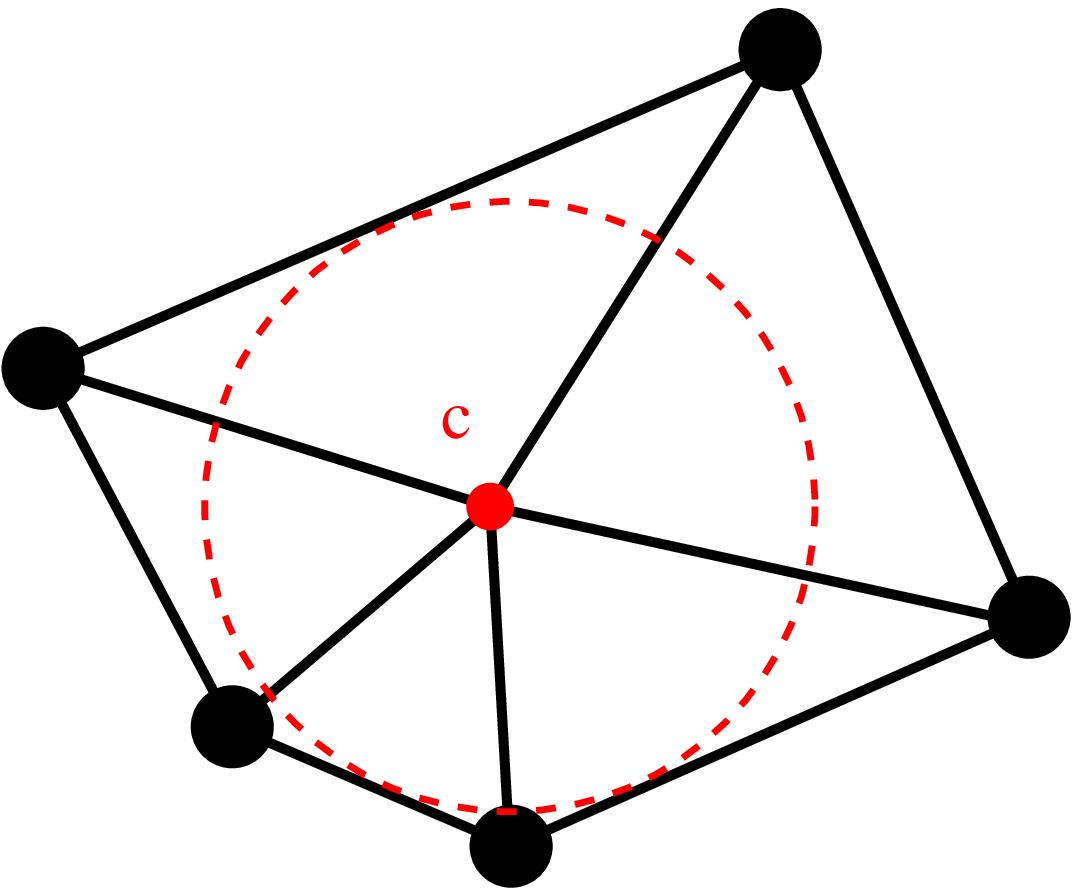} &
\includegraphics[height=.35\linewidth]{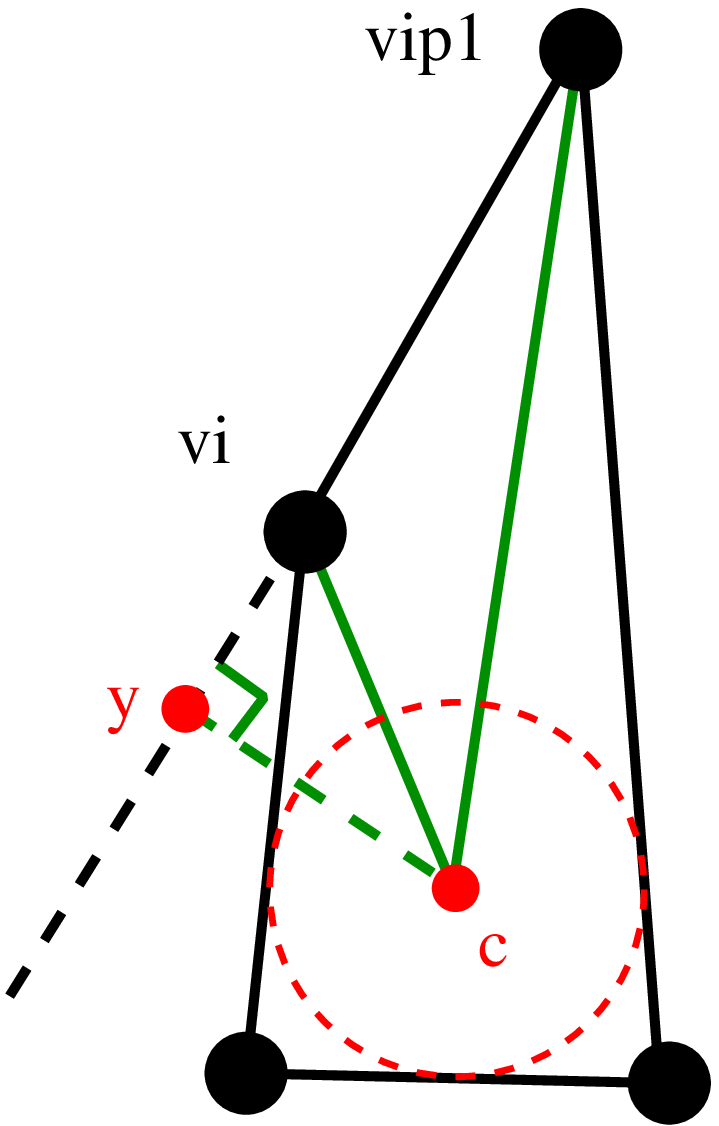} &
\text{ }\quad\text{ } &
\includegraphics[height=.35\linewidth]{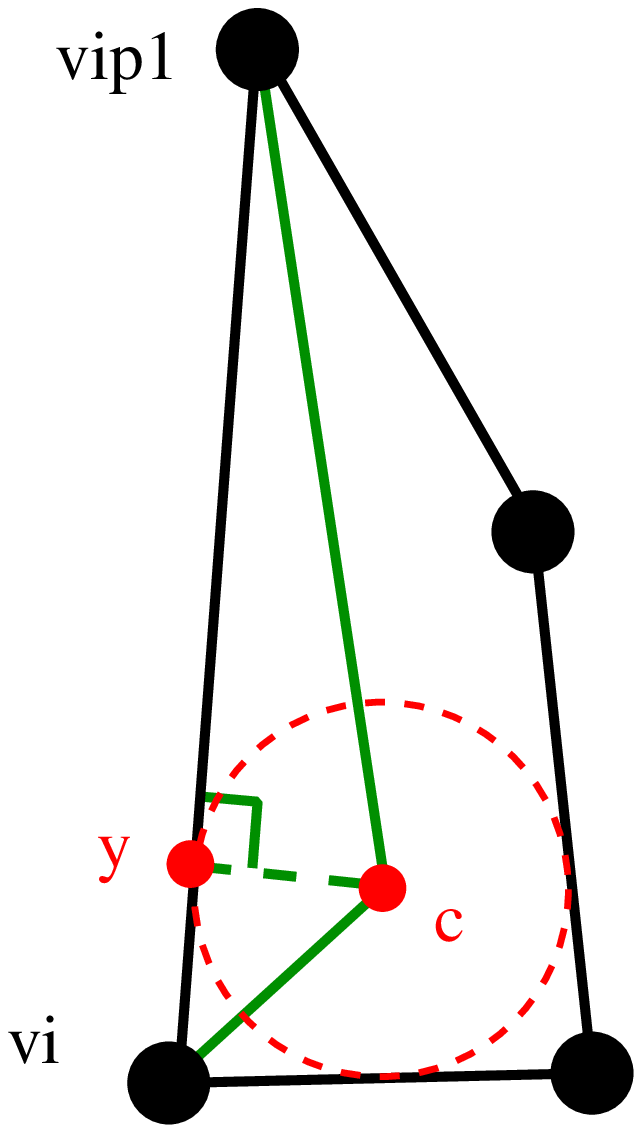} \\
\\
\textbf{(a)} &
\textbf{(b)} &&
\textbf{(c)} 
\end{array}\]
\end{center} 
\caption{\textbf{(a)} Triangulation used in the analysis of Harmonic coordinates. \textbf{(b)} Notation for proof of the bound for $\angle \bc\bv_i\bv_{i+1}$ in a case where it is $>\pi/2$. \textbf{(c)} Notation for proof of the bound for $\angle \bv_i\bc\bv_{i+1}$ in a case where it is $>\pi/2$.}
\label{fig:opt-tri}
\end{figure}

Recalling the notation from Figure~\ref{fig:notation}, let $T$ be the triangulation of $\Omega$ formed by connecting $\bc$ to each of the $\bv_i$; see Figure~\ref{fig:opt-tri}a.  

\begin{proposition}\label{prop:nlatriangulation}
Under G\ref{g:ratio} all angles of all triangles of $T$ are less than $\pi - \arcsin (1/\gamma^*)$.  
\end{proposition}
\begin{proof}
Consider the triangle with vertices $\bc$, $\bv_i$ and $\bv_{i+1}$.  Without loss of generality, assume that $|\bc - \bv_i| < |\bc- \bv_{i+1}|$.  First we bound $\angle\bc \bv_{i+1} \bv_i$.  By the law of sines,
\begin{equation}
\label{eq:lawsinarg}
\frac{\sin(\angle\bc \bv_{i+1} \bv_i)}{\sin(\angle\bc \bv_i \bv_{i+1})}=\frac{|\bc - \bv_i|}{|\bc - \bv_{i+1}|}<1.
\end{equation}
If $\angle\bc \bv_i \bv_{i+1}>\pi/2$ then $\angle\bc \bv_{i+1} \bv_i<\pi/2$.  Otherwise, (\ref{eq:lawsinarg}) implies $\angle\bc\bv_{i+1}\bv_i<\pi/2$.

To bound angle $\angle \bc\bv_i\bv_{i+1}$, it suffices to consider the case when $\angle \bc\bv_i\bv_{i+1}>\pi/2$, as shown in Figure~\ref{fig:opt-tri}b.  Define $\by$ to be the point on the line through $\bv_i$ and $\bv_{i+1}$ which forms a right triangle with $\bv_i$ and $\bc$.  Since $\angle \bc\bv_i\bv_{i+1}>\pi/2$, $\by$ is exterior to $\Omega$, as shown.  Observe that
\[\frac{|\bc-\bv_i|}{|\bc-\by|}<\frac{|\bc-\bv_{i+1}|}{|\bc-\by|}<\frac{\diam(\Omega)}{\rho(\Omega)}=\gamma<\gamma^*.\]
Since $\sin(\pi-\angle \bc\bv_i\bv_{i+1})=\frac {|\bc-\by|}{|\bc-\bv_i|}$, the result follows.

For the final case, it suffices to assume $\angle \bv_i\bc\bv_{i+1} > \pi/2$, as shown in Figure~\ref{fig:opt-tri}c.  Define $\by$ in the same way, but note that in this case $\by$ is between $\bv_i$ and $\bv_{i+1}$, as shown.  Similarly, $\frac{|\bc-\bv_{i+1}|}{|\bc-\by|}<\gamma^*$, implying $\angle\bv_i\bv_{i+1}\bc >\arcsin(1/\gamma^*)$.  Since $\angle\bv_i\bc\bv_{i+1}<\pi-\angle\bv_i\bv_{i+1}\bc$, the result follows.\qed


\end{proof}

\begin{lemma}\label{lem:optimal}
Under G\ref{g:ratio} the operator $\Iopt$ satisfies the $H^1$ interpolant estimate (\ref{eq:iuh1uh2}).
\end{lemma}
\begin{proof}
Since the differential equation (\ref{eq:optpde}) is linear, $\Iopt u$ is the solution to the differential equation,
\begin{equation}
\label{eq:ioptpde}
\ds\left\{\begin{array}{rcll}
\Delta\left(\Iopt u \right) & = & 0, & \text{on $\Omega$} \\
\Iopt_i u & = & g_u, & \text{on $\p\Omega$}
\end{array}\right.
\end{equation}
where $g_u$ is the piecewise linear function which equals $u$ at the vertices of $\Omega$.  Following the standard approach for handling nonhomogeneous boundary data we divide: $\Iopt u = u_{\rm hom} + u_{\rm non}$ where $u_{\rm non}\in H^1(\Omega)$ is some function satisfying the boundary condition (i.e., $u_{\rm non} = g_u$ on $\p\Omega$) and $u_{\rm hom}$ solves,
\begin{equation}
\label{eq:uhom}
\ds\left\{\begin{array}{rcll}
\Delta u_{\rm hom} & = & -\Delta u_{\rm non}, & \text{on $\Omega$} \\
u_{\rm hom} & = & 0, & \text{on $\p\Omega$}.
\end{array}\right.
\end{equation}
Specifically we select $u_{\rm non}$ to be the standard Lagrange interpolant of $u$ over triangulation $T$ (described earlier).
Since Proposition~\ref{prop:nlatriangulation} guarantees that no large angles exist in the triangulation, the standard interpolation error estimate holds,
\begin{equation}
\label{eq:nlainterpolation}
\hpn{u - u_{\rm non}}{1}{\Omega} \leq C_{BA} \, \hpsn{u}{2}{\Omega}
\end{equation}
where $C_{BA}$ only depends upon the aspect ratio bound $\gamma^*$ since $\diam(\Omega) = 1$.  
The triangle inequality then implies that $\hpn{u_{\rm non}}{1}{\Omega} \leq \max(1,C_{BA}) \hpn{u}{2}{\Omega}$.

Next a common energy estimate (see \cite{Ev98}) for (\ref{eq:uhom}) implies that $\hpsn{u_{\rm hom}}{1}{\Omega} \leq \hpsn{u_{\rm non}}{1}{\Omega}$.  The Poincar\'e inequality (see \cite{Le09}) ensures that $\lpn{u_{\rm hom}}{2}{\Omega} \leq C_P \hpsn{u}{1}{\Omega}$ where $C_P$ only depends on the diameter of $\Omega$ which we have fixed to be $1$.   The argument is completed by combining the previous estimates:
\begin{align*}
\hpn{\Iopt u}{1}{\Omega} & \leq \hpn{u_{\rm hom}}{1}{\Omega} + \hpn{u_{\rm non}}{1}{\Omega} \\
 & \leq (1 + C_P) \hpsn{u_{\rm hom}}{1}{\Omega} + \hpn{u_{\rm non}}{1}{\Omega}\\
& \leq (1 + C_P)\hpn{u_{\rm non}}{1}{\Omega} \leq (1 + C_P)\max(1,C_{BA}) \hpn{u}{2}{\Omega}.  
\qed 
\end{align*}
\end{proof}

\subsection{Estimate Requirements for Wachspress Coordinates}
\label{sec:wachpress}

\begin{figure}[ht]
\begin{center}
\psfrag{v1}{$\bv_1$}
\[\begin{array}{ccc}
\includegraphics[width=.3\linewidth]{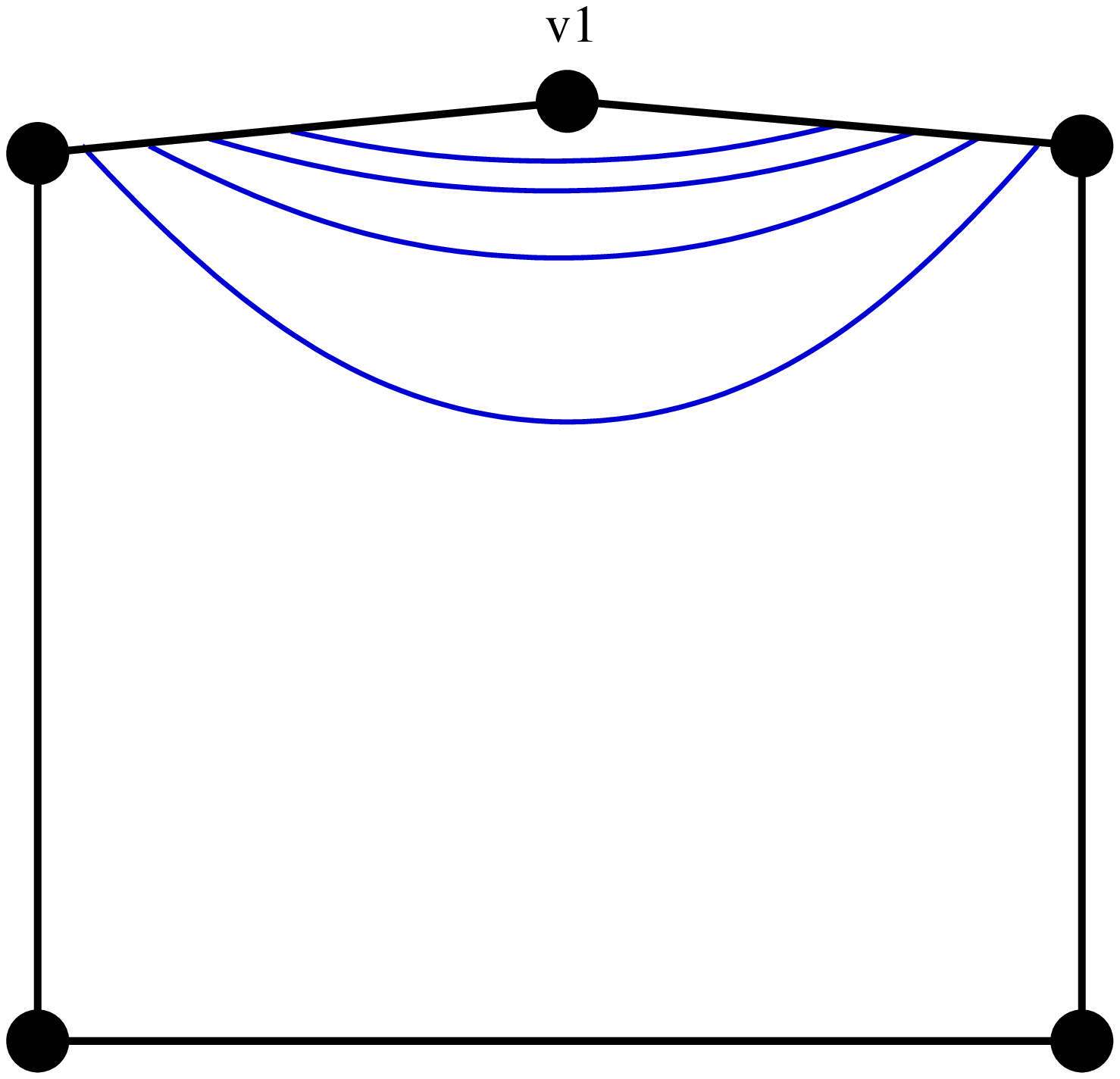} &
\includegraphics[width=.3\linewidth]{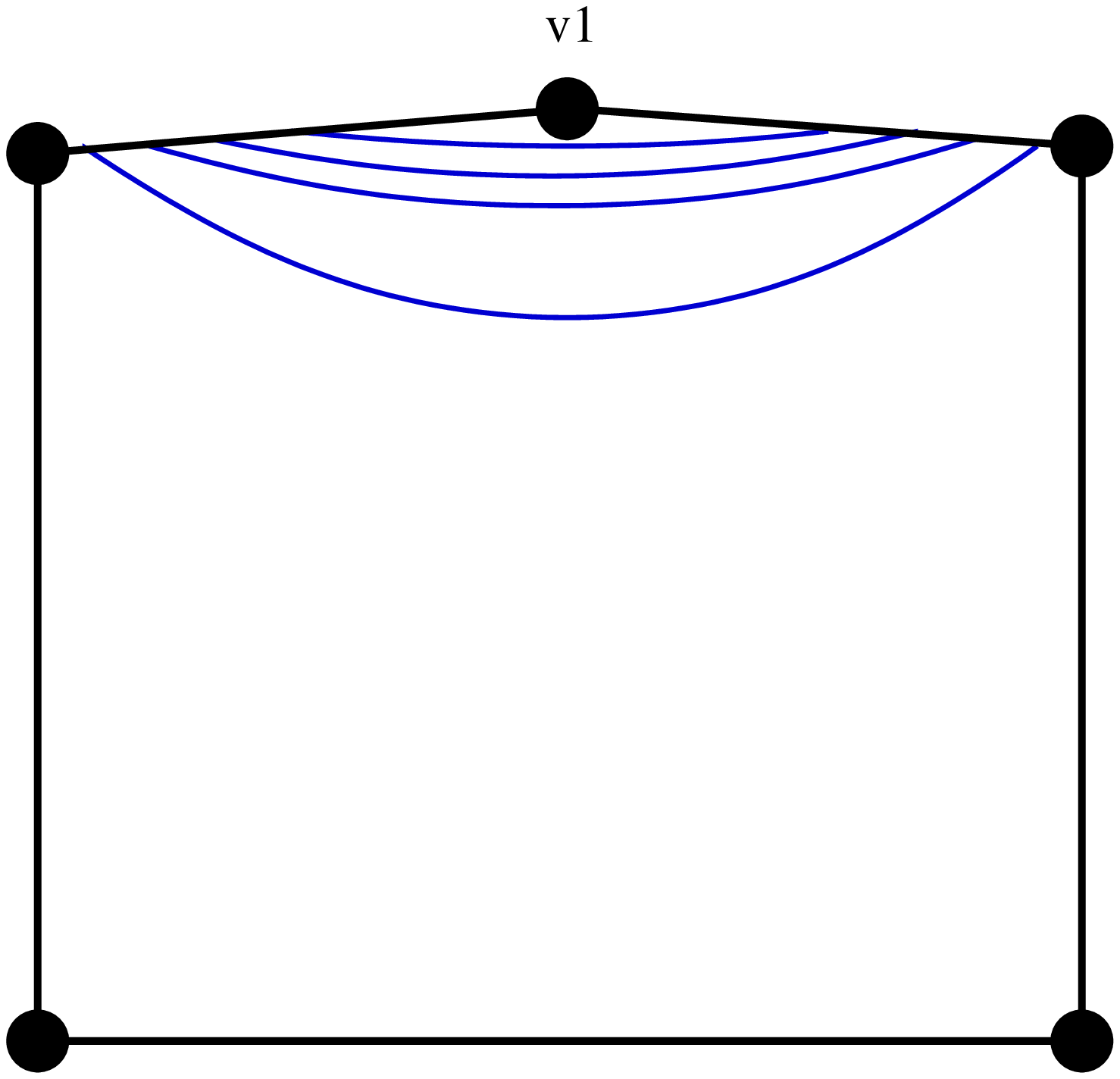} &
\includegraphics[width=.3\linewidth]{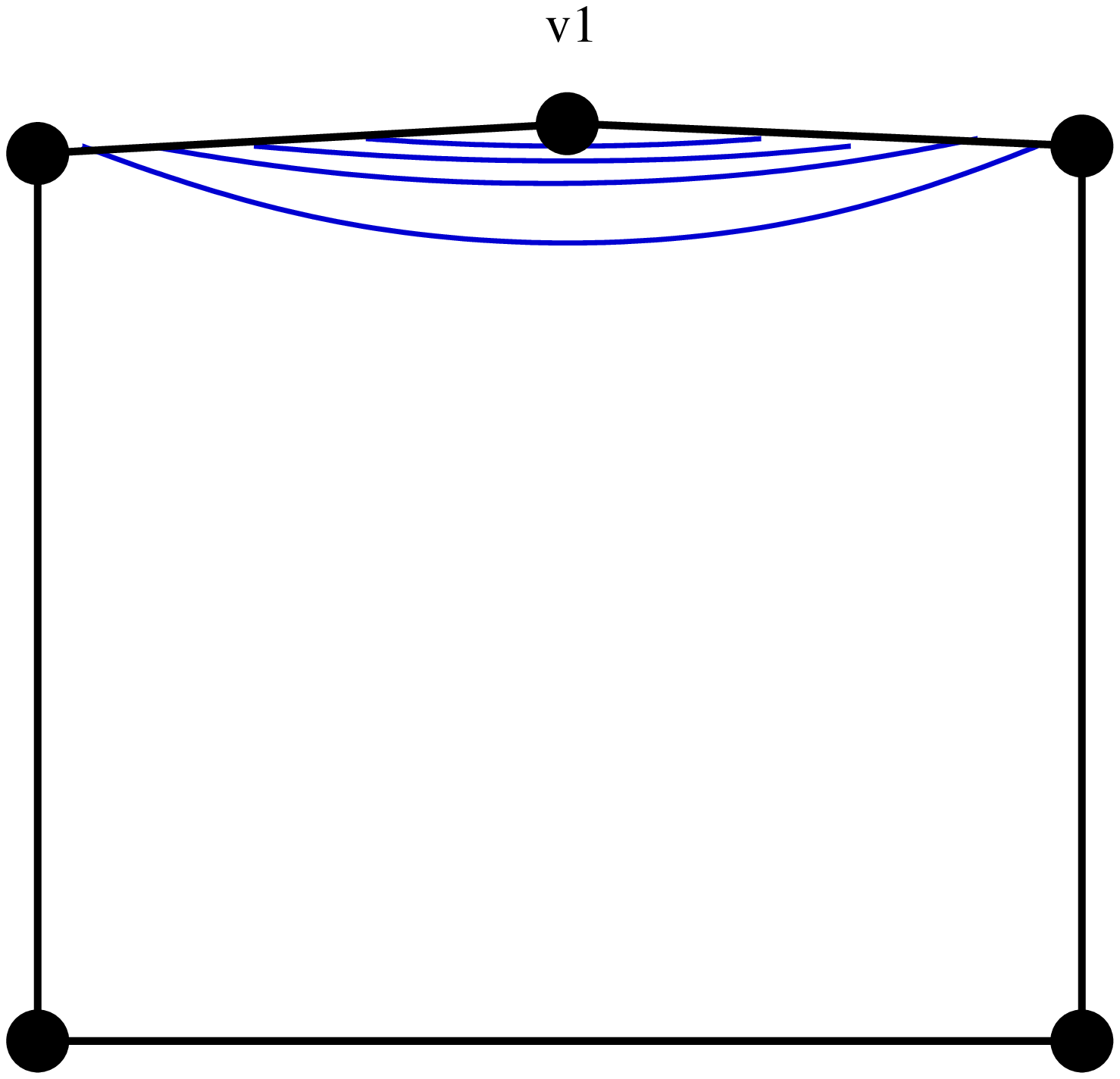}
\end{array}\]
\end{center} 
\caption{Example showing the necessity of condition G\ref{g:maxangle} for attaining the optimal convergence estimate (\ref{eq:hconv}) with the Wachspress coordinates.  As the shape approaches a square, the level sets of $\lwach_1$ collect at the top edge, causing a steep gradient and thus preventing a bound on the $H^1$ norm of the error.  The figures from left to right correspond to $\eps$ values of 0.1, 0.05, and 0.025. }
\label{fig:wach-cex}
\end{figure}

Unlike the harmonic coordinate functions, the Wachspress coordinates can produce unsatisfactory interpolants unless additional geometric conditions are imposed.  We present a simple counterexample (observed qualitatively in \cite{FHK2006} and in Figure~\ref{fig:wach-cex}) to show what can go wrong.  

Let $\pent$ be the pentagon defined by the vertices
\[
\bv_1=(0,1+\eps),\quad
\bv_2=(1,1),\quad
\bv_3=(1,-1),\quad
\bv_4=(-1,-1),\quad
\bv_5=(-1,1),
\]
with $\eps>0$.  As $\eps\raw 0$, $\pent$ approaches a square so G1 is not violated.  
Consider the interpolant of $u(\bx)=1-x_1^2$ where $\bx =(x_1,x_2)$.  Observe that $u$ has value 1 at $\bv_1$ and value $0$ at the other vertices of $\pent$.  Hence
\[\Iwach u=\sum_{i=1}^5 u(\bv_i)\lwach_i=\lwach_1\]
Using the fact that $\p u/\p y = 0$, we write
\begin{align*}
\hpn{u - \Iwach u}{1}{\pent}^2 &= \hpn{u - \lwach_1}{1}{\pent}^2 \\
 & = \int_\pent |u - \lwach_1|^2+\left|\frac{\p (u - \lwach_1)}{\p x}\right|^2+\left|\frac{\p \lwach_1}{\p y}\right|^2.
\end{align*}
The last term in this sum blows up as $\epsilon \rightarrow 0$.  
\begin{lemma}
\label{lem:wachcex}
$\ds\lim_{\eps\raw 0}\int_\pent\left|\frac{\p \lwach_1}{\p y}\right|^2=\infty$.
\end{lemma}
The proof of the lemma is given in the appendix.  As a corollary, we observe that $\hpn{u - Iu}{1}{\pent}$ cannot be bounded independently of $\eps$.  Since $\hpn{u}{2}{\pent}$ is finite, this means the optimal convergence estimate (\ref{eq:hconv}) cannot hold without additional geometric criteria on the domain $\Omega$.  This establishes the necessity of a maximum interior angle bound on vertices if Wachspress coordinates are used.

Under the three geometric restrictions G\ref{g:ratio}, G\ref{g:minedge}, and G\ref{g:maxangle}, (\ref{eq:hconv}) does hold which will be shown in Lemma~\ref{lm:wach}.  We begin with some preliminary estimates.

\begin{proposition}\label{pr:gradarea}
For all $\bx \in \Omega$, $\vsn{\nabla A_i(\bx)} \leq \frac{1}{2}$.
\end{proposition}
\begin{proof}
In \cite[Equation (17)]{FK10} it is shown that the $\vsn{\nabla A_i(\bx)}=\frac{1}{2}\vsn{\bv_i - \bv_{i+1}}$.  Since $\diam (\Omega) = 1$ the result follows.  
\qed
\end{proof}

Next we show that the triangular areas $B_i$ are uniformly bounded from below given our geometric assumptions.  

\begin{figure}[ht]
\begin{center}
\psfrag{vi}{$\bv_i$}
\psfrag{vj}{$\bv_j$}
\psfrag{vjp1}{$\bv_{j+1}$}
\psfrag{vip1}{$\bv_{i+1}$}
\psfrag{vim1}{$\bv_{i-1}$}
\psfrag{dstar}{$d_*$}
\psfrag{x}{$\bx$}
\psfrag{Ajx}{$A_j(\bx)$}
\[
\includegraphics[width=1.8in]{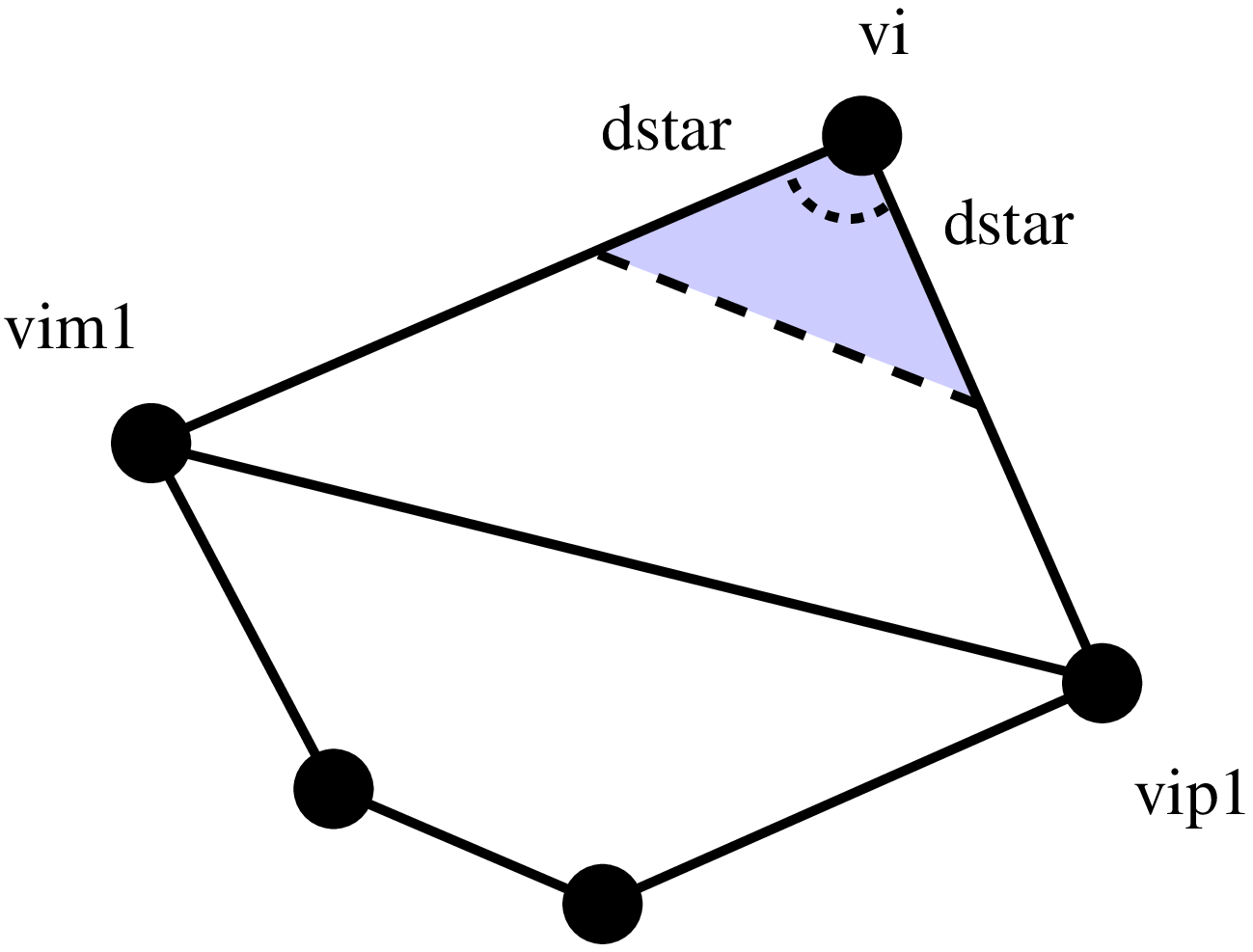}\quad\text{ }\quad
\includegraphics[width=1.5in]{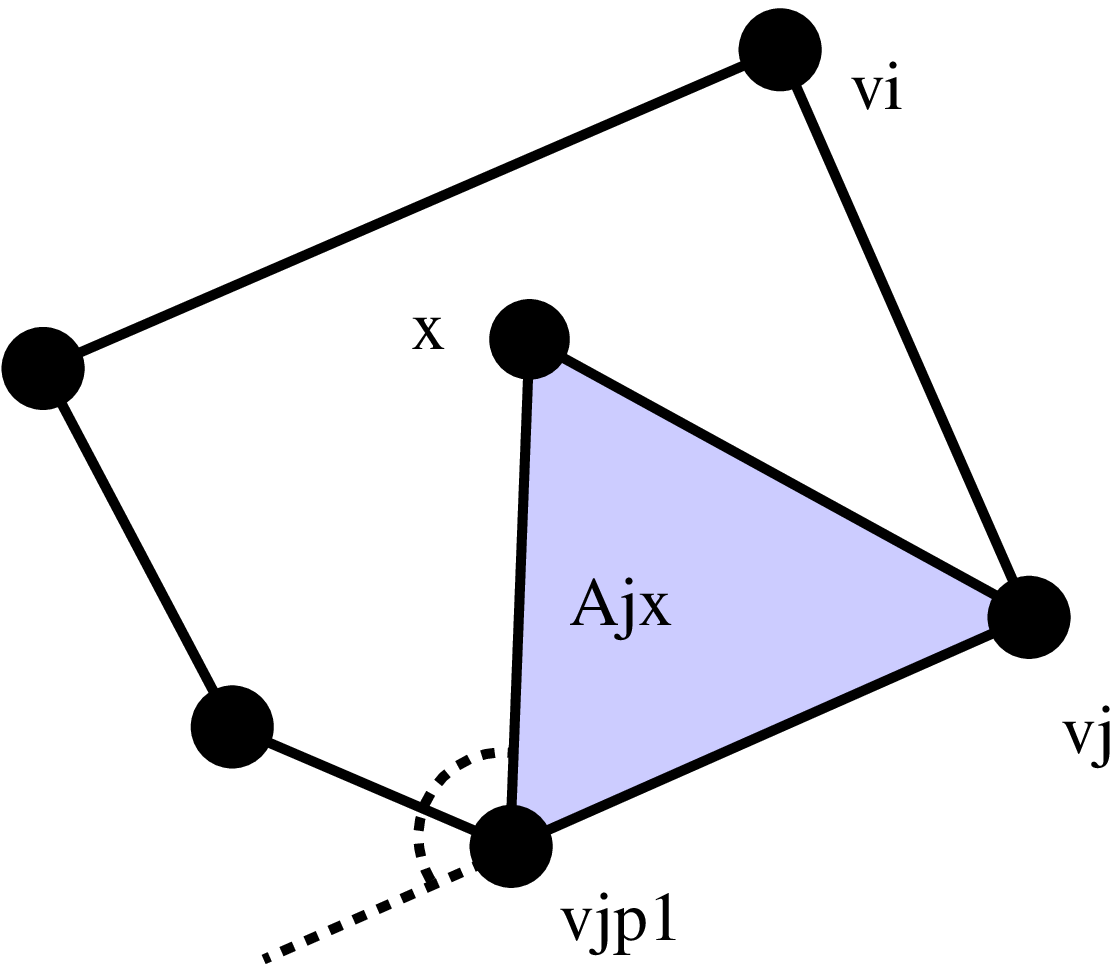}
\]
\end{center} 
\caption{\textbf{Left:} Justification of claim that $B_i > (d_*)^2 \sin (\beta_i/2) \cos(\beta_i/2)$ in the proof of Proposition~\ref{pr:carealower}.  The shaded triangle is isosceles with angle $\beta_i$ and two side lengths equal to $d_*$ as indicated.  Computing the area of this triangle using the dashed edge as the base yields the estimate.  \textbf{Right:} Justification of claim that $A_j(\bx) > (d_*)^2\sin(\pi-\beta^*)/4$ in the proof of Proposition~\ref{pr:wachweightlower}.  The indicated angle is at least $\pi-\beta^*$ by G\ref{g:maxangle} and $\vsn{\bv_{j}-\bv_{j+1}}>d_*$.  Computing the area of the triangle using edge $\bv_{j}\bv_{j+1}$ as the base yields the estimate.}
\label{fig:ajxlower}
\end{figure}

\begin{proposition}\label{pr:carealower}
Under G\ref{g:ratio}, G\ref{g:minedge}, and G\ref{g:maxangle}, there exists $B_*$ such that $B_i >B_*$.  
\end{proposition}
\begin{proof}
By G\ref{g:minedge}, the area of the isosceles triangle with equal sides of length $d_*$ meeting with angle $\beta_i$ at $\bv_i$ is a lower bound for $B_i$, as shown in Figure~\ref{fig:ajxlower} (left).  More precisely, $B_i > (d_*)^2 \sin (\beta_i/2) \cos(\beta_i/2)$.  G\ref{g:maxangle} implies that $\cos(\beta_i/2) > \cos(\beta^*/2)$.  G\ref{g:minangle} (which follows from G\ref{g:ratio} by Proposition \ref{prop:Grelation}) implies that $\sin (\beta_i/2) > \sin(\beta_*/2)$. Thus $B_i > B_* := (d_*)^2 \sin (\beta_*/2) \cos(\beta^*/2)$.  
\qed
\end{proof}

Proposition~\ref{pr:carealower} can be extended to guarantee a uniform lower bound on the sum of the Wachspress weight functions.  

\begin{proposition}\label{pr:wachweightlower}
Under G\ref{g:ratio}, G\ref{g:minedge}, and G\ref{g:maxangle}, there exists $w_*$ such that for all $\bx\in \Omega$, $$\sum_k \wwach_k(\bx) > w_*.$$
\end{proposition}
\begin{proof}
Let $\bv_i$ be the nearest vertex to $\bx$, breaking any tie arbitrarily.  We will produce a lower bound on $w_i(\bx)$.  Let $j \notin \{i-1,i\}$.  G\ref{g:minedge} implies that $\vsn{\bx-\bv_j} > d_*/2$ and $\vsn{\bx-\bv_{j+1}} > d_*/2$.  G\ref{g:maxangle} implies that $\angle \bx\bv_j\bv_{j+1} < \beta^*$ and $\angle \bx\bv_{j+1}\bv_j < \beta^*$.  It follows that $A_j(\bx) > (d_*)^2\sin(\pi-\beta^*)/4$ (see Figure~\ref{fig:ajxlower} (right)).  We now use  Proposition~\ref{pr:carealower} and property G\ref{g:maxdegree} (which follows from either G\ref{g:minedge} or G\ref{g:maxangle} by Proposition~\ref{prop:Grelation}) to conclude that
\[
\sum_k \wwach_k(\bx) > \wwach_i(\bx) = B_i \prod_{j\not=i,i-1}A_j(\bx) \geq B_* \left[(d_*)^2\sin(\pi-\beta^*)/4\right]^{n^*-2}.
\]
\qed
\end{proof}

\begin{lemma}\label{lm:wach}
Under G1, G2, and G3, (\ref{eq:basisbound}) holds for the Wachspress coordinates.
\end{lemma}
\begin{proof}

The gradient of $\lwach_i(\bx)$ can be bounded using Proposition \ref{pr:wachweightlower}:
\begin{align}
\vsn{\nabla \lwach_i(\bx)} & \leq \frac{\vsn{\nabla \wwach_i(\bx)}}{\sum_j \wwach_j(\bx)} + \frac{\wwach_i(\bx)\sum_k \vsn{\nabla \wwach_k(\bx)}}{\left(\sum_j \wwach_j(\bx)\right)^2}\notag\\
 & \leq \frac{\vsn{\nabla \wwach_i(\bx)} + \sum_k \vsn{\nabla \wwach_k(\bx)}}{\sum_j \wwach_j(\bx)}
 \leq \frac{2 \sum_k \vsn{\nabla \wwach_k(\bx)}}{w_*}.\label{eq:wachgradbound}
\end{align}
Recalling Proposition~\ref{prop:convexfacts}, $\sum_{j=1}^{n} A_j(\bx) < \pi/4$ and $B_i<\pi/4$.
Using Proposition~\ref{pr:gradarea} and the arithmetic mean-geometric mean inequality, we derive
\begin{align}
\vsn{\nabla \wwach_i(\bx)} &= \vsn{\sum_{j\neq i-1,i} B_i \nabla A_j(\bx) \prod_{k\neq i-1,i,j} A_k(\bx)}\notag \\
 & \leq \sum_{j\neq i-1,i}\left[ B_i \vsn{\nabla A_j(\bx)} \prod_{k\neq i-1,i,j} A_k(\bx)\right]\notag \\
 & \leq \sum_{j\neq i-1,i}\frac{\pi}{8}\left[ \frac{\sum_{k\neq i-1,i,j} A_k(\bx)}{n-3} \right]^{n-3} 
 \leq \sum_{j\neq i-1,i}\frac{\pi}{8}\left[ \frac{\pi}{4(n-3)} \right]^{n-3} \notag \\
 & = \frac{\pi}{8}(n-2)\left[ \frac{\pi}{4(n-3)}\right]^{n-3}. \label{eq:wachgradweight}
\end{align}
By induction, one can show that $n(n-2)\left[ \frac{\pi}{4(n-3)}\right]^{n-3}  \leq 2\pi$ for $n\geq 4$.  Using this, we substitute (\ref{eq:wachgradweight}) into (\ref{eq:wachgradbound}) to get
\[\vsn{\nabla \lwach_i(\bx)}\leq \frac 2{w_*}\sum_k|\nabla\wwach_k(\bx)|\leq\frac 2{w_*}n\frac{\pi}{8}(n-2)\left[ \frac{\pi}{4(n-3)}\right]^{n-3}\leq \frac{\pi}{4w_*}2\pi =\frac{\pi^2}{2w_*}. \]
Since $|\Omega|<\pi/4$ by Proposition \ref{prop:convexfacts}, we thus have a uniform bound
\[\hpn{\lwach_i}{1}{\Omega} \leq \sqrt{\left(1 + \frac{\pi^4}{4w_*^2}\right)\frac{\pi}{4}}.\qed \] 
\end{proof}

\subsection{Estimate Requirements for Sibson Coordinates}
\label{sec:sibson}

The interpolation estimate for Sibson coordinates is computed using a very similar approach to that of the previous section on Wachspress coordinates.  However in this case the geometric condition G\ref{g:maxangle} is not necessary.  We begin with a technical property of domains satisfying conditions G\ref{g:ratio} and G\ref{g:minedge}.  

\begin{proposition}\label{prop:ballinvoronoi}
Under G\ref{g:ratio} and G\ref{g:minedge}, there exists $h_* > 0$ such that for all $\bx\in \Omega$, $B(\bx,h_*)$ does not intersect any three edges or any two non-adjacent edges of $\Omega$.  
\end{proposition}
\begin{proof}
Let $\bx \in \Omega$, $h\in (0,d_*/2)$, and suppose that two disjoint edges of $\Omega$, $e_i$ and $e_j$, intersect $B(\bx, h)$.  Let $L_i$ and $L_j$ be the lines containing $e_i$ and $e_j$ and let $\theta$ be the angle between these lines; see Figure~\ref{fig:ballinset}.  We first consider the case where $L_i$ and $L_j$ are not parallel and define $\bz = L_i \cap L_j$.

\begin{figure}
\centering
\psfrag{x}{$\bx$}
\psfrag{z}{$\bz$}
\psfrag{ei}{$e_i$}
\psfrag{ej}{$e_j$}
\psfrag{vi}{$\bv_i$}
\psfrag{vj}{$\bv_j$}
\psfrag{Li}{$L_i$}
\psfrag{Lj}{$L_j$}
\psfrag{h}{$h$}
\psfrag{W}{$W$}
\psfrag{theta}{$\theta$}
\psfrag{omega}{$\Omega$}
\includegraphics[width=.75\linewidth]{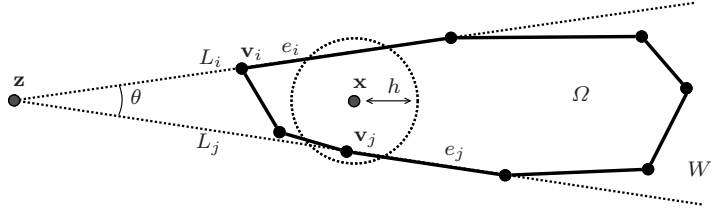}
\caption{Notation for proof of Proposition~\ref{prop:ballinvoronoi}.}
\label{fig:ballinset}
\end{figure}

Let $\bv_i$ and $\bv_j$ be the endpoints of $e_i$ and $e_j$ nearest to $\bz$.  Since $h < d_*/2$ both $\bv_i$ and $\bv_j$ cannot live in $B(\bx, h)$; without loss of generality assume that $\bv_i\notin B(\bx, h)$.  

Since $\dist(\bv_j,L_i) < 2h$, 
\begin{equation}\label{eq:thetabound}
\sin\theta < 2h/\vsn{\bz - \bv_j}.
\end{equation}

Let $W$ be the sector between $L_i$ and $L_j$ containing $x$.  Now $\Omega \subset B(\bv_j,1)\cap W \subset B(\bz, 1 + \vsn{\bz-\bv_j})\cap W$.  
It follows that $\rho(\Omega) \leq (1 + \vsn{\bv_j - \bz})\sin\theta $.  
Using (\ref{eq:thetabound}) and G\ref{g:ratio}, 
\[
\frac{1}{\gamma^*} \leq \frac{2h}{\vsn{\bz - \bv_j}}(1 + \vsn{\bz - \bv_j}) \leq 2h\left(\frac{1}{d_*} + 1 \right)
\]
where the final inequality holds because by G\ref{g:minedge}
$\vsn{\bz - \bv_j} \geq \vsn{\bv_i - \bv_j} \geq d_*$.
Thus
\begin{equation}\label{eq:hlwrbd}
h > \frac{d_*}{2\gamma^*(1+d_*)}.
\end{equation}
Estimate (\ref{eq:hlwrbd}) holds in the limiting case: when $L_i$ and $L_j$ are parallel.  In this case $\Omega$ must be contained in a strip of width $2h$ which for small $h$ violates the aspect ratio condition.  

The triangle is the only polygon with three or more pairwise non-adjacent edges.  
So it remains to find a suitable $h_*$ so that $B(\bx,h_*)$ does not intersect all three edges of the triangle.  
For a triangle, $\rho(\Omega)$ is the radius of the smallest circle touching all three edges.  
Since under G\ref{g:ratio} $\rho(\Omega) \geq 1/\gamma^*$, $B(\bx,\frac{1}{2\gamma^*})$ intersects at most two edges.  
Thus $h_* = \frac{d_*}{2\gamma^*(1+d_*)}$ is sufficiently small to satisfy the proposition in all cases.  
\qed
\end{proof}

Proposition~\ref{prop:ballinvoronoi} is a useful tool for proving a lower bound on $D(\bx)$, the area of the Voronoi cell of $\bx$ intersected with $\Omega$.

\begin{proposition}\label{prop:minvoronoiarea}
Under G\ref{g:ratio} and G\ref{g:minedge}, there exists $D_* >0$ such that $D(\bx) > D_*$.  
\end{proposition}
\begin{proof}
Let $h_*$ be the constant in Proposition~\ref{prop:ballinvoronoi}.  We consider two cases, based on whether the point $\bx$ is near any vertex of $\Omega$, as shown in Figure~\ref{fig:sibsproof} (left). 

\noindent \underline{Case 1}: There exists $\bv_i$ such that $\bx \in B(\bv_i,h_*/2)$.

\begin{figure}
\centering
\psfrag{x}{$\bx$}
\psfrag{vi}{$\bv_i$}
\psfrag{hstar}{$h_*$}
\psfrag{Bi}{$\beta_i$}
\[\begin{array}{ccc}
\includegraphics[width=.3\linewidth]{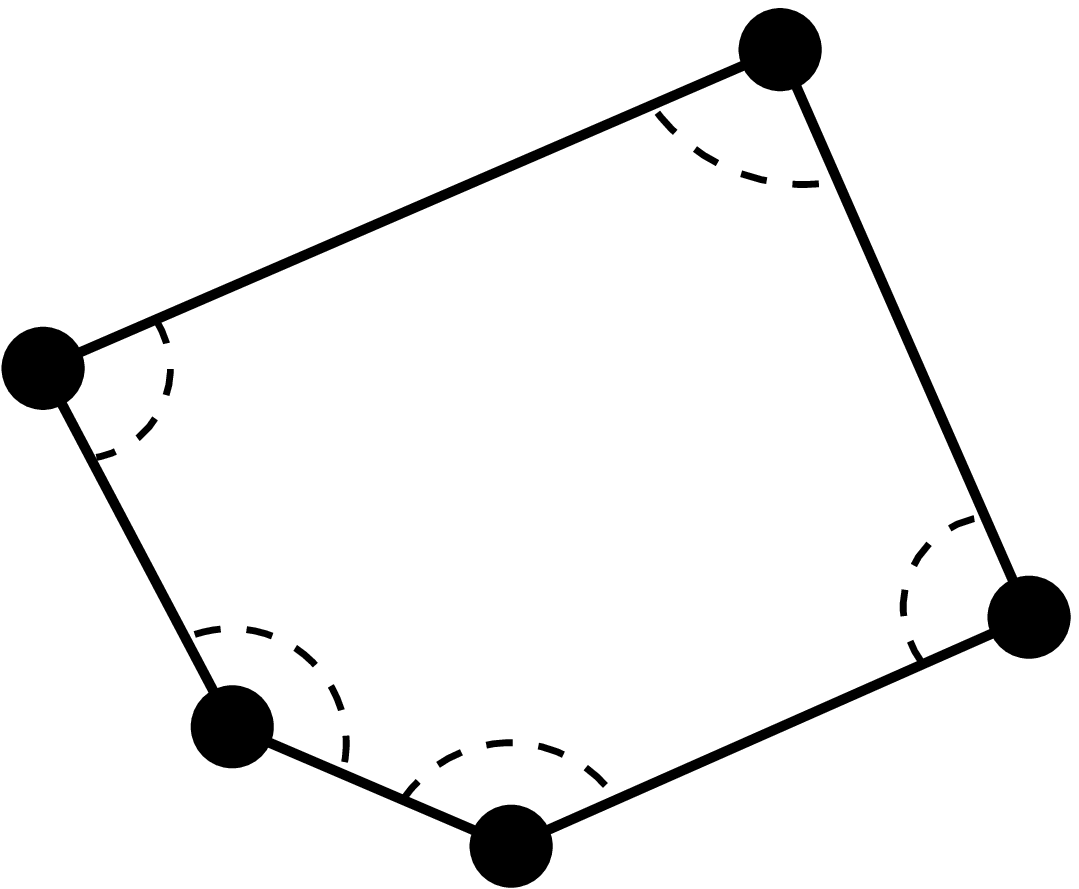} &
\quad\text{  }\quad &
\includegraphics[width=.23\linewidth]{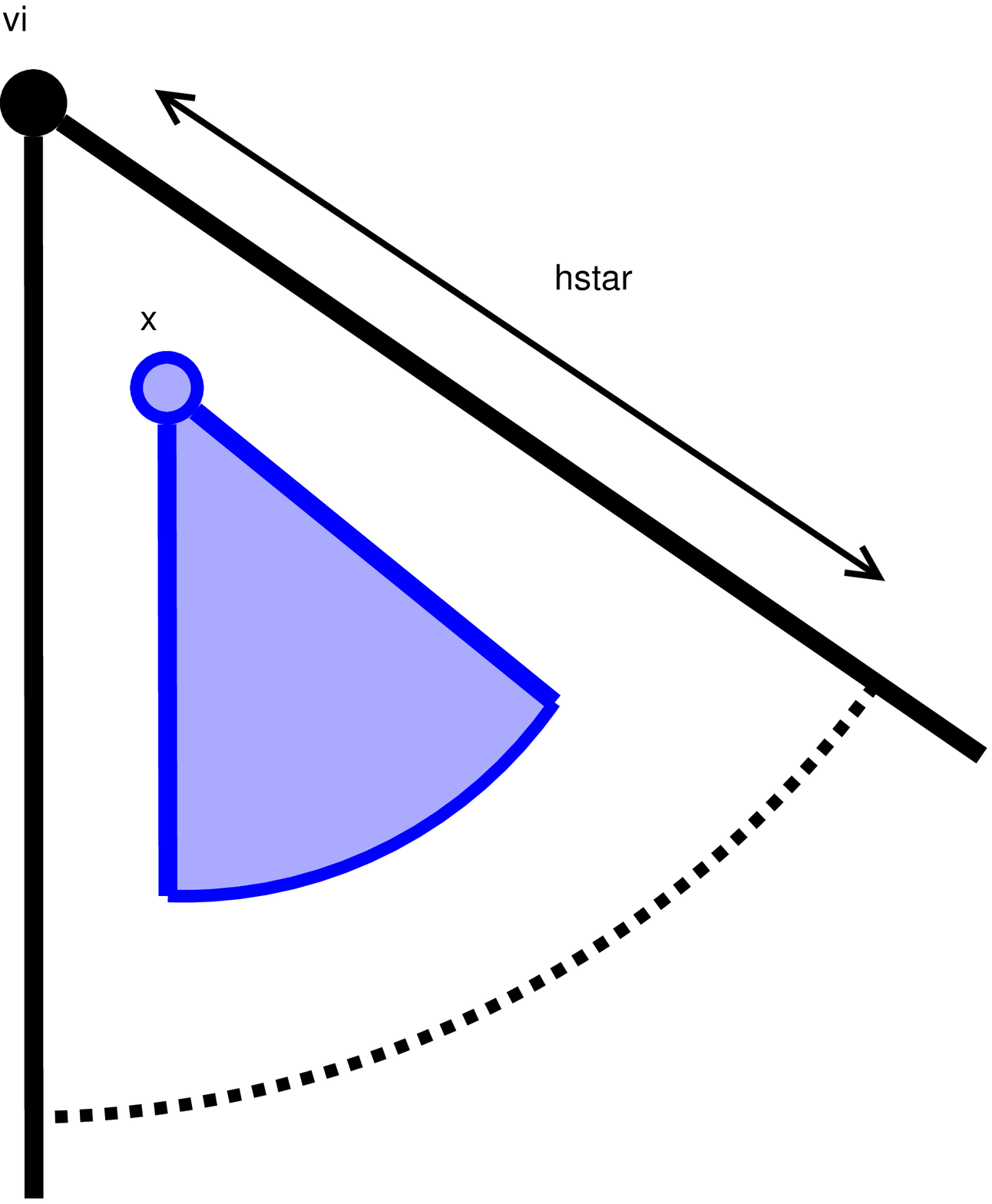}
\end{array}\]
\caption{The proof of Proposition~\ref{prop:minvoronoiarea} has two cases based on whether $\bx$ is within $h_*/2$ of some $\bv_i$ or not.  When $\bx$ is within $h_*/2$ of $\bv_i$, the shaded sector shown on the right  is contained in $V_{P'}(\bx)\cap\Omega$.}
\label{fig:sibsproof}
\end{figure}

Consider the sector of $B(\bx,h_*/2)$ specified by segments which are parallel to the edges of $\Omega$ containing $\bv_i$, as shown in Figure~\ref{fig:sibsproof} (right).  This sector must be contained in $\Omega$ by Proposition~\ref{prop:ballinvoronoi} and in the Voronoi cell of $\bx$ by choice of $h_*<d_{*}$.  Thus by G\ref{g:minangle} (using Proposition~\ref{prop:Grelation}(\ref{gr:getminangle})) $D(\bx) \geq \beta_* h_*^2/8$.

\noindent \underline{Case 2}: For all $\bv_i$, $\bx \notin B(\bv_i,h_*/2)$.

In this case, $B(\bx, h_*/4) \cap \Omega \subset V_{P'}(\bx)$.  If $B(\bx, h_*/4)$ intersects zero or one boundary edge of $\Omega$, then $D(\bx) \geq \pi h_*^2/32$.  Otherwise $B(\bx, h_*/4)$ intersects two adjacent boundary edges. By G\ref{g:minangle}, $D(\bx) \geq \beta_* h_*^2/32$.
\qed
\end{proof}

General formulas for the gradient of the area of a Voronoi cell are well-known and can be used to bound the gradients of $D(\bx)$ and $D(\bx)\cap C_i$.

\begin{proposition}\label{prop:voronoigrad}
$|\nabla D(\bx)| \leq \pi$ and $|\nabla (D(\bx)\cap C_i)| \leq 1$.
\end{proposition}
\begin{proof}
The gradient of the area of a Voronoi region is known to be
\[
\nabla D(\bx) = \sum_{j=1}^{n} \frac{\bv_j - \bx}{\vsn{\bv_j - \bx}} F_j,
\]
where $F_j$ is the length of the segment separating the Voronoi cells of $\bx$ and $\bv_j$~\cite{OA91,OBSC00}.  Then applying Proposition~\ref{prop:convexfacts} gives 
\[
\vsn{\nabla D(\bx)} \leq \sum_{i=1}^{n} F_i \leq |\partial\Omega|\leq  \pi.
\]
Similarly, 
\[
\nabla (D(\bx)\cap C_i) = \frac{\bv_i - \bx}{\vsn{\bv_i - \bx}} F_i,
\]
and since $F_i \leq \diam(\Omega)$, $|\nabla (D(\bx)\cap C_i)| \leq 1$.
\qed
\end{proof}

Propositions~\ref{prop:minvoronoiarea} and \ref{prop:voronoigrad} give estimates for the key terms needed in proving (\ref{eq:basisbound}) for the Sibson coordinates $\lsib$.  

\begin{lemma}\label{lm:sibs}
Under G1 and G2, (\ref{eq:basisbound}) holds for the Sibson coordinates.
\end{lemma}
\begin{proof}
$\vsn{\nabla \lsib_i}$ is estimated by applying Propositions~\ref{prop:minvoronoiarea} and \ref{prop:voronoigrad}:
\begin{align*}
\vsn{\nabla \lsib_i} &\leq \frac{|\nabla (D(\bx)\cap C_i)|}{D(\bx)} + \frac{(D(\bx)\cap C_i) \vsn{\nabla D(\bx)} }{D(\bx)^2}
 \leq \frac{|\nabla (D(\bx)\cap C_i)| + |\nabla D(\bx)|}{D(\bx)} \\
&\leq \frac{1 + \pi}{D_*}.  
\end{align*}
Integrating this estimate completes the result.  
\qed
\end{proof}
\begin{corollary}
By Lemma \ref{lem:basisbd}, the $H^1$ interpolant estimate (\ref{eq:iuh1uh2}) holds for the Sibson coordinates.
\end{corollary}

\section{Final Remarks}
\label{sec:conc}

Geometric requirements needed to ensure optimal interpolation error estimates are necessary for guaranteeing the compatibility of polygonal meshes with generalized barycentric interpolation schemes in finite element methods.  Moreover, the identification of necessary and unnecessary geometric restrictions provides a tool for comparing various approaches to barycentric interpolation.  Specifically we have demonstrated the necessity of a maximum interior angle restriction for Wachspress coordinates, which was empirically observed in \cite{FHK2006}, and shown that this restriction is unneeded when using Sibson coordinates.  

Table \ref{tab:conditions} provides a guideline for how to choose barycentric basis functions given geometric criteria or, conversely, which geometric criteria should be guaranteed given a choice of basis functions.  While utilized throughout our analysis, the aspect ratio requirement G\ref{g:ratio} can likely be substantially weakened.  Due to a dependence on specific affine transformations, such techniques on triangular domains \cite{BA76,Ja76} (i.e., methods for proving error estimates under the maximum angle condition rather than the minimum angle condition) cannot be naturally extended to polygonal domains.  Although aimed at a slightly different setting that we have analyzed, challenges in identifying sharp geometric restrictions are apparent from the numerous studies on quadrilateral elements, e.g., \cite{Ja77,ZV95,AD01,MNS08}.  A satisfactory generalization of the maximum angle condition to arbitrary polygons is a subject of further investigation.

This paper emphasizes three specific barycentric coordinates (harmonic, Wachspress, and Sibson) but several others have been introduced in the literature.  Maximum entropy \cite{Su04}, metric \cite{MLBD02}, and discrete harmonic \cite{PP93} coordinates can all be studied either by specific analysis or generalizing the arguments given here to wider classes of functions. The mean value coordinates defined by Floater~\cite{F2003} are of particular interest in this regard as they are defined by an explicit formula and appear to not require a maximum angle condition.  The formal analysis of these functions, however, is not trivial.  Additional generalizations could be considered by dropping certain restrictions on the coordinates, such as non-negativity, or the mesh elements, such as convexity.  Working with non-convex elements, however, would require some non-obvious generalization of the geometric restrictions G1-G5.  



\appendix

\section{Proof of Lemma \ref{lem:wachcex}}
\noindent
\begin{proof}
The explicit formula for the Wachspress weight associated to $\bv_1$ is
\[\wwach_1(\bx)= B_1A_2A_3A_4 =\eps(1-x)(1+x)(1+y) \]
where $\bx=(x,y)$ is an arbitrary point inside $\pent$.  The other weights can be computed similarly, yielding the coordinate function
\[\lwach_1= \frac{w_{i}(\bv)}{\sum_{j=1}^5w_{j}(\bv)} = \frac{\eps(1-x)(1+x)(1+y)}{\eps^2 (1 - x^2) + 4\eps + 2(1 - y)}.\]
The $y$ partial derivative term is computed to be
\[\frac{\p \lwach_1}{\p y} = \frac{4\eps(1-x^2) + 4\eps^2(1-x^2) + \eps^3(1 - x^2)^2}{(\eps^2 (1 - x^2) + 4\eps + 2(1 - y))^2}\]
Define the subregion $\pentp\subset\pent$ by
\[\pentp=\left\{(x,y)\in\pent:\frac 14\leq x\leq\frac 34,\quad 1\leq y\leq 1+\eps\right\}\]
Observe that $\frac{7}{16}\leq 1 - x^2 \leq \frac{15}{16}$ on $\pentp$.  Fix $0<\eps<1$.  We bound the numerator by
\[
4\eps(1-x^2) + 4\eps^2(1-x^2) + \eps^3(1 - x^2)^2 > 4\eps\cdot \frac{7}{16} + 4 \eps^2 \cdot \frac{7}{16} + \eps^3\cdot \frac{49}{256} > \frac{7}{4}\eps.
\]
Since $|y-1|<\eps$ on $\pentp$, we can bound the denominator by
\[
|\eps^2 (1 - x^2) + 4\eps + 2(1 - y)| \leq |\eps^2 (1 - x^2)| + |4\eps| + |2(1 - y)| \leq \eps^2 + 4\eps + 2\eps \leq 7\eps.
\]
Putting these results together, we have that
\[\left|\frac{\p \lwach_1}{\p y}\right|>\frac{\frac{7}{4}\eps}{49\eps^2} = \frac{1}{28\eps }>0.\]
Let $C=\frac {1}{28}$ for ease of notation.  Since $|\pentp|>\frac 18 \eps$,
\[\lim_{\eps\raw 0}\int_\pent\left|\frac{\p \lwach_1}{\p y}\right|^2\geq \lim_{\eps\raw 0}\int_{\pentp}\left|\frac{\p \lwach_1}{\p y}\right|^2>\lim_{\eps\raw 0}\int_{\pentp}\frac{C^2}{\eps^2}=C^2\lim_{\eps\raw 0}\frac{|\pentp|}{\eps^2}>\frac{C^2}{8}\lim_{\eps\raw 0}\frac 1\eps=\infty,\]
thereby proving the lemma.
\qed
\end{proof}

\bibliographystyle{spmpsci}      
\bibliography{journal}

\end{document}